\documentclass[3p,12pt]{elsarticle}
\usepackage{stmaryrd}
\usepackage{amsfonts}
\usepackage{amsmath,amssymb}
\usepackage{mathrsfs}
\usepackage{amstext}
\usepackage{subeqnarray}
\usepackage{txfonts,graphics,url,cases,fancyvrb}
\usepackage{multirow}

\newtheorem{example}{Example}[section]

\newtheorem{definition}{Definition}[section]

\newtheorem{remark}{Remark}[section]
\newenvironment{proof}{\par{\bf Proof.}}{\qquad\endproof}
\def\endproof{\vbox{\hrule height0.6pt\hbox{%
 \vrule height1.3ex width0.6pt\hskip0.8ex
 \vrule width0.6pt}\hrule height0.6pt}}

\newcommand{\dx}[1][x]{\,{\rm d}#1}

\begin{document}

\begin{frontmatter}

\title{Fractional differentiation matrices with applications}

\author{Fanhai Zeng, Changpin Li}


\begin{abstract}
In this paper, the fractional differential matrices based on the
Jacobi-Gauss points are derived with respect to the Caputo and Riemann-Liouville
fractional derivative operators.  The spectral radii of the fractional differential matrices
are investigated numerically. The spectral collocation schemes are illustrated
to solve the fractional ordinary differential equations and fractional partial
differential equations.
Numerical examples are also presented to illustrate the effectiveness
of the derived methods, which show better performances over some
existing methods.
\end{abstract}

\end{frontmatter}


\section{Introduction}\label{sec:1}
Differential matrices are useful and easily implemented in the simulation
of the classical differential equations \cite{Canuto2006,ShenTangWang2011}.
This paper aims to develop the fractional
differential matrices to approximate the fractional integral and derivative operators
with applications for solving the fractional differential equations.

Fractional calculus (including the fractional integral  and the
fractional derivative) has   become a hot topic recently for its
wide applications in many areas of science and engineering, see for example
\cite{Cortes2007,Kulish2002,LiuYangTurner2011,Oldham1974,Podlubny1999,RossikhinShitikova2010,Samko93,
YuLiu2008,ZhuangGuLiu2011}.

Unlike the classical derivative operator, the fractional derivative operators
are nonlocal with weakly singularity, which are more complicated
for theoretical analysis and numerical simulation. Up to now, there have
been some numerical methods to discretize the fractional integral and derivative
operators, for instance
\cite{CaoXu2013,DiethelmFordFreed2004,DiethelmFord2005,GaoLiao2011,GaoSun2014,KumarAgrawal2006,
LakestaniDehghan2012,LiChenYe2011b,LiZhaoChen2011,Lubich86,LynchCarreras2003,
Oldham1974,Murio2006,Sousa2012,SugiuraHasegawa2009,Odibat2006,Odibat2009}.
In \cite{Langlands2005,SunWu2006}, the L1 method was proposed to discretize the
Riemann-Liouville and Caputo derivative  operators.
Tian et al. \cite{TianZhaouDeng2012} proposed the weighted formula based on the shifted
Gr\"{u}nwald-Letnikov formula to approximate the Riemann--Liouville derivative operator with
second-order accuracy. \c{C}elik and Duman \cite{CelikDuman2012} proposed the fractional central difference
method to discretize Riesz fractional derivative operator with convergence of order 2.
The operational matrices based on the explicit forms of the Legendre, Chebyshev and Jacobi polynomials
were proposed to discretize the Caputo derivative operator in
\cite{DohaBhrawyEzz-Eldien2011a,DohaBhrawyEzz-Eldien2011,DohaBhrawyEzz-Eldien2012,
EsmaeiliShamsi2011,SaadatmandiDehghan2010}. Tian and Deng \cite{TianDeng2013}
proposed a method to approximate the
Caputo fractional derivative operator with the fractional differential matrices obtained,
but their method seems unstable when performing on the common computers with double precision.
Xu and Hesthaven \cite{XuHesthaven2014} also proposed the fractional differential matrix
to approximate the Caputo derivative operator, but the derivation of the fractional differential
matrix involves in calculating the inverse matrix. Recently, a multi-domain spectral method based on the
multi-domain fractional differential matrix for
time-fractional differential equation was developed in \cite{ChenXuHesthaven14}.

In \cite{LiZengLiu2012}, the effective recurrence formulas were developed to approximate
the fractional integral and the left Caputo fractional derivative of the Legendre,
Chebyshev, and Jacobi polynomials, and the corresponding operational  matrices are obtained
such that the fractional integral and derivative of a given function at one collocation
point can be calculated with $O(N)$ operations.  In this paper, we choose the collocation
points as the Jacobi-Gauss types, such that the fractional differentiation
of a given function $u(x),x\in [-1,1]$ at the collocation points $x_j\,(j=0,1,...,N)$
are approximated by the matrix-vector product,
i.e., $A\mathbf{u}$, $A\in \mathbb{R}^{(N+1)\times (N+1)},\mathbf{u}=(u(x_0),u(x_1),...,u(x_N))^T$.
We call such a type of matrix $A$ the fractional differential matrix. The spectral radius  $\rho(A)$
of the fractional differential matrix $A$  is numerically
investigated, which shows the behavior as $\rho(A){\,\leq\,}C_0N^{2\alpha}$, where $C_0$ is independent
of $N$, $\alpha$ is the order of the corresponding fractional derivative operator.
The spectral collocation schemes are illustrated to solve the fractional ordinary differential
equations and the  fractional partial differential equations. Numerical experiments display good
satisfactory results, and the comparison between other methods are made to show better performances
of the present methods.

The remainder of this paper is organized as follows. In Section 2,
we introduce several definitions of fractional calculus  and the
Legendre, Chebyshev and Jacobi polynomials. The fractional differential matrices
with respect to the Caputo and Riemann--Liouville derivative operators are derived
in Section 3. The spectral collocation methods for solving the fractional differential equations
are illustrated in Section 4. Numerical examples are
presented in Section 5, and the conclusion is included in the last section.

\section{Preliminaries}
In this section,  we introduce the definitions of the fractional
calculus. Then we introduce the Legendre, Chebyshev and Jacobi
polynomials, which will be used later on.

\begin{definition}\label{DefInt}
The left and right fractional integrals (or the left and right Riemann--Liouville integrals)
with  order $\alpha>0$ of the given function $f(t)$ are
defined as
\begin{equation} \label{RLInt1}
D^{-\alpha}_{a,t}f(t)=\frac{1}{\Gamma(\alpha)}\int_a^t(t-s)^{\alpha-1}f(s)\dx[s]
\end{equation}
and
\begin{equation} \label{RLInt1-2}
D^{-\alpha}_{t,b}f(t)=\frac{1}{\Gamma(\alpha)}\int_t^b(s-t)^{\alpha-1}f(s)\dx[s]
\end{equation}
respectively, where $\Gamma(\cdot)$ is the Euler's gamma function.
\end{definition}

There exist several kinds of fractional derivatives, which will be introduced in the following.
\begin{definition}\label{DefRL}
The  left and right Riemann-Liouville fractional  derivatives with  order $\alpha>0$ of the given function
$f(t),t\,\in\,(a,b)$ are defined  as
\begin{equation} \label{rlfd1}
 \begin{aligned}
{}_CD_{a,t}^{\alpha}f(t)=\frac{\mathrm{d}^n}{\mathrm{d}t^n}D^{-(n-\alpha)}_{a,t}f(t)
=\frac{1}{\Gamma(n-\alpha)}\frac{\mathrm{d}^n}{\mathrm{d}t^n}\int_a^t(t-s)^{n-\alpha-1}f(s)\dx[s]
 \end{aligned}
\end{equation}
and
\begin{equation} \label{rlfd2}
 \begin{aligned}
{}_CD_{t,a}^{\alpha}f(t)=(-1)^n\frac{\mathrm{d}^n}{\mathrm{d}t^n}D^{-(n-\alpha)}_{t,a}f(t)
=\frac{(-1)^n}{\Gamma(n-\alpha)}\frac{\mathrm{d}^n}{\mathrm{d}t^n}\int_t^b(s-t)^{n-\alpha-1}f(s)\dx[s]
 \end{aligned}
\end{equation}
respectively, where $n$ is a positive integer and $n-1<\alpha{\,\leq\,} n$.
\end{definition}

\begin{definition}\label{DefCapto}
The  left and right Caputo fractional   derivatives with  order $\alpha>0$ of the given function
$f(t),t\,\in\,(a,b)$ are defined  as
\begin{equation} \label{cptfd1}
 \begin{aligned}
{}_CD_{a,t}^{\alpha}f(t)=D^{-(n-\alpha)}_{a,t}\Big[f^{(n)}(t)\Big]
=\frac{1}{\Gamma(n-\alpha)}\int_a^t(t-s)^{n-\alpha-1}f^{(n)}(s)\dx[s]
 \end{aligned}
\end{equation}
and
\begin{equation} \label{cptfd2}
 \begin{aligned}
{}_CD_{t,a}^{\alpha}f(t)=(-1)^nD^{-(n-\alpha)}_{t,a}\Big[f^{(n)}(t)\Big]
=\frac{(-1)^n}{\Gamma(n-\alpha)}\int_t^b(s-t)^{n-\alpha-1}f^{(n)}(s)\dx[s]
 \end{aligned}
\end{equation}
respectively, where $n$ is a positive integer and $n-1<\alpha{\,\leq\,} n$.
\end{definition}

\begin{definition}\label{DefRZ}
The  Riesz fractional derivative and Riesz-Caputo fractional derivative  with  order $\alpha>0$
of the given function $f(t),t\,\in\,(a,b)$ are defined  as
\begin{equation} \label{rzfrac}
{}_{RZ}D_t^{\alpha}f(t)=c_{\alpha}\Big({}_{RL}D_{a,t}^{\alpha}f(t)+{}_{RL}D_{t,b}^{\alpha}f(t)\Big),
\end{equation}
and
\begin{equation} \label{rzcaptofrac}
{}_{RC}D_t^{\alpha}f(t)=c_{\alpha}\Big({}_{C}D_{a,t}^{\alpha}f(t)+{}_{C}D_{t,b}^{\alpha}f(t)\Big),
\end{equation}
respectively, where $c_{\alpha}=-\frac{1}{2\cos(\frac{\pi\alpha}{2})}$.
\end{definition}

%

Next, we introduce the Jacobi polynomials. The Jacobi polynomials $\{P^{a,b}_j(x)\},\,a,b>-1,x\in[-1,1]$ are
given by the following three-term recurrence relation
\cite{ShenTangWang2011}
\begin{equation} \label{jacobi}
 \begin{aligned}
&P^{a,b}_0(x)=1,{\quad}P^{a,b}_1(x)=\frac{1}{2}(a+b+2)x+\frac{1}{2}(a-b),\\
&P^{a,b}_{j+1}(x)=(A_j^{a,b}x-B_j^{a,b})P^{a,b}_{j}(x)-C_j^{a,b}P^{a,b}_{j-1}(x),\quad
{n\geq 1},
 \end{aligned}
\end{equation}
where
\begin{equation} \label{jacobi1}
 \begin{aligned}
&A_j^{a,b}=\frac{(2j+a+b+1)(2j+a+b+2)}{2(j+1)(j+a+b+1)},\\
&B_j^{a,b}=\frac{(b^2-a^2)(2j+a+b+1)}{2(j+1)(j+a+b+1)(2j+a+b)},\\
&C_j^{a,b}=\frac{(j+a)(j+b)(2j+a+b+2)}{(j+1)(j+a+b+1)(2j+a+b)}.
 \end{aligned}
\end{equation}

Next, several properties of the Jacobi polynomials will be stated. Let
$\omega^{a,b}(x)=(1-x)^a(1+x)^b.$
Then, one has
\begin{equation} \label{jacobi1-2}
\int_{-1}^{1}P^{a,b}_{m}(x)P^{a,b}_{n}(x)\omega^{a,b}(x)\dx
=\left\{\begin{aligned}
&0,{\quad}m\neq n\\
&\gamma^{a,b}_n,{\quad}m=n,
 \end{aligned}\right.
\end{equation}
where
\begin{equation} \label{jacobi1-3}
\gamma^{a,b}_n=\frac{2^{a+b+1}\Gamma(n+a+1)\Gamma(n+b+1)}{(2n+a+b+1)n!\Gamma(n+a+b+1)}.
\end{equation}

Some other properties of the Jacobi polynomials that
will be used in the present paper  are presented below:
\begin{equation} \label{jacobi2}
 \begin{aligned}
&P^{a,b}_{j}(1)=\binom{j+a}{j}=\frac{\Gamma(j+a+1)}{j!\Gamma(a+1)},{\quad}
P^{a,b}_{j}(-1)=(-1)^j\frac{\Gamma(j+b+1)}{j!\Gamma(b+1)}.
 \end{aligned}
\end{equation}

\begin{equation} \label{jacobi3}
 \begin{aligned}
\frac{{\rm
d}^m}{\dx^m}P^{a,b}_{j}(x)=d_{j,m}^{a,b}P^{a+m,b+m}_{j-m}(x),\quad
j{\,\geq\,}m,m\in\mathbb{N},
 \end{aligned}
\end{equation}
where
\begin{equation} \label{jacobi3.2}
d_{j,m}^{a,b}=\frac{\Gamma(j+m+a+b+1)}{2^m\Gamma(j+a+b+2)}.
\end{equation}

\begin{equation} \label{jacobi4}
 \begin{aligned}
P^{a,b}_{j}(x)=\widehat{A}_j^{a,b}\frac{{\rm
d}}{\dx}P^{a,b}_{j-1}(x)+\widehat{B}_j^{a,b}\frac{{\rm
d}}{\dx}P^{a,b}_{j}(x)+ \widehat{C}_j^{a,b}\frac{{\rm
d}}{\dx}P^{a,b}_{j+1}(x),\quad {j\geq 1},
 \end{aligned}
\end{equation}
in which
\begin{equation} \label{jacobi5}
 \begin{aligned}
&\widehat{A}_j^{a,b}=\frac{-2(j+a)(j+b)}{(j+a+b)(2j+a+b)(2j+a+b+1)},\\
&\widehat{B}_j^{a,b}=\frac{2(a-b)}{(2j+a+b)(2j+a+b+2)},\\
&\widehat{C}_j^{a,b}=\frac{2(j+a+b+1)}{(2j+a+b+1)(2j+a+b+2)}.
 \end{aligned}
\end{equation}
If $j=1$, then $\widehat{A}_1^{a,b}=0$ in \eqref{jacobi5}.

If $a=b=0$, then the recurrence formula \eqref{jacobi} is reduced to the Legendre polynomials as
\begin{equation} \label{lgd}
 \begin{aligned}
&L_0(x)=1,{\quad}L_1(x)=x,{\quad}
L_{j+1}=\frac{2j+1}{j+1}xL_j(x)-\frac{j}{j+1}L_{j-1}(x),\quad
{j\geq 1}.
 \end{aligned}
\end{equation}
If $a=b=-1/2$ in \eqref{jacobi}, then
$P^{-\frac{1}{2},-\frac{1}{2}}_{j}(x)=\frac{\Gamma(j+1/2)}{j!\sqrt{\pi}}T_j(x)$,
where $T_j(x)$ is the Chebyshev polynomial that can be defined as
\begin{equation} \label{cheb}
T_0(x)=1,{\quad}T_1(x)=x, {\quad}T_{j+1}(x)=2xT_j(x)-T_{j-1}(x),\quad {j\geq 1}.
\end{equation}

\section{Derivations of the fractional differential matrices}
Denote by
$$\widehat{P}_{L,j}^{a,b,\alpha}(x)=\frac{1}{\Gamma(\alpha)}\int_{-1}^x(x-s)^{\alpha-1}P^{a,b}_{j}(s)\dx[s].$$
Using the recurrence formulae \eqref{jacobi}--\eqref{jacobi1}, the
properties \eqref{jacobi2} and \eqref{jacobi4}, we have
\begin{equation} \label{eq:fint1}
\left\{\begin{aligned}
&\widehat{P}_{L,0}^{a,b,\alpha}(x)=\frac{(x+1)^{\alpha}}{\Gamma(\alpha+1)},\\
&\widehat{P}_{L,1}^{a,b,\alpha}(x)=\frac{a+b+2}{2}\left(\frac{x(x+1)^{\alpha}}{\Gamma(\alpha+1)}
-\frac{\alpha(x+1)^{\alpha+1}}{\Gamma(\alpha+2)}\right)+\frac{a-b}{2}\widehat{P}_{L,0}^{a,b,\alpha}(x),\\
&\widehat{P}_{L,j+1}^{a,b,\alpha}(x)=\frac{A_j^{a,b}x-B_j^{a,b}-{\alpha}A_j^{a,b}\widehat{B}_j^{a,b}}
{1+{\alpha}A_j^{a,b}\widehat{C}_j^{a,b}}\widehat{P}_{L,j}^{a,b,\alpha}(x)
-\frac{C_j^{a,b}+{\alpha}A_j^{a,b}\widehat{A}_j^{a,b}}{1+{\alpha}A_j^{a,b}\widehat{C}_j^{a,b}}\widehat{P}_{L,j-1}^{a,b,\alpha}(x)\\
&{\quad}+\frac{{\alpha}A_j^{a,b}\left(\widehat{A}_j^{a,b}P^{a,b}_{j-1}(-1)+\widehat{B}_j^{a,b}P^{a,b}_{j}(-1)
+\widehat{C}_j^{a,b}P^{a,b}_{j+1}(-1)\right)}
{\Gamma(\alpha+1)\left(1+{\alpha}A_j^{a,b}\widehat{C}_j^{a,b}\right)}
(x+1)^{\alpha},{\quad\quad}j{\,\geq\,}1.
\end{aligned}\right.
\end{equation}

The above recurrence formula was first obtained in \cite{LiZengLiu2012}. Denote by
$$\widehat{P}_{R,j}^{a,b,\alpha}(x)=\frac{1}{\Gamma(\alpha)}\int_{x}^1(s-x)^{\alpha-1}P^{a,b}_{j}(s)\dx[s].$$
For $j{\,\geq\,}1$, we can obtain from \eqref{jacobi} that
\begin{equation} \label{eq:fint1-2}
 \begin{aligned}
\widehat{P}_{R,j+1}^{a,b,\alpha}(x)
=&\frac{1}{\Gamma(\alpha)}\int_{x}^1(s-x)^{\alpha-1}P^{a,b}_{j+1}(s)\dx[s]\\
=&\frac{1}{\Gamma(\alpha)}\int_{x}^1(s-x)^{\alpha-1}
\left[(A_j^{a,b}s-B_j^{a,b})P^{a,b}_{j}(s)-C_j^{a,b}P^{a,b}_{j-1}(s)\right]\dx[s]\\
=&(A_j^{a,b}x-B_j^{a,b})\widehat{P}_{R,j}^{a,b,\alpha}(x) -C_j^{a,b}\widehat{P}_{R,j-1}^{a,b,\alpha}(x)
+ \frac{A_j^{a,b}}{\Gamma(\alpha)}\int_{x}^1(s-x)^{\alpha}P^{a,b}_{j}(s)\dx[s].
 \end{aligned}
\end{equation}
Using \eqref{jacobi4} yields
\begin{equation} \label{eq:fint1-3}
 \begin{aligned}
\widehat{P}_{R,j+1}^{a,b,\alpha}(x)
=&(A_j^{a,b}x-B_j^{a,b})\widehat{P}_{R,j}^{a,b,\alpha}(x) -C_j^{a,b}\widehat{P}_{R,j-1}^{a,b,\alpha}(x)\\
&+\frac{A_j^{a,b}}{\Gamma(\alpha)}\int_{x}^1(s-x)^{\alpha}\,\mathrm{d}\left[\widehat{A}_j^{a,b}P^{a,b}_{j-1}(s)
+\widehat{B}_j^{a,b}P^{a,b}_{j}(s)+ \widehat{C}_j^{a,b}P^{a,b}_{j+1}(s)\right]\\
=&(A_j^{a,b}x-B_j^{a,b})\widehat{P}_{R,j}^{a,b,\alpha}(x) -C_j^{a,b}\widehat{P}_{R,j-1}^{a,b,\alpha}(x)\\
&+\frac{A_j^{a,b}}{\Gamma(\alpha)}(s-x)^{\alpha}\left[\widehat{A}_j^{a,b}P^{a,b}_{j-1}(s)
+\widehat{B}_j^{a,b}P^{a,b}_{j}(s)+ \widehat{C}_j^{a,b}P^{a,b}_{j+1}(s)\right]\bigg|_{x}^{1}\\
&-\frac{{\alpha}A_j^{a,b}}{\Gamma(\alpha)}\int_{x}^1(s-x)^{\alpha-1}\left[\widehat{A}_j^{a,b}P^{a,b}_{j-1}(s)
+\widehat{B}_j^{a,b}P^{a,b}_{j}(s)+ \widehat{C}_j^{a,b}P^{a,b}_{j+1}(s)\right]\dx[s]\\
=&(A_j^{a,b}x-B_j^{a,b})\widehat{P}_{R,j}^{a,b,\alpha}(x) -C_j^{a,b}\widehat{P}_{R,j-1}^{a,b,\alpha}(x)\\
&-{\alpha}A_j^{a,b}\left(\widehat{A}_j^{a,b}\widehat{P}^{a,b}_{R,j-1}(s)
+\widehat{B}_j^{a,b}\widehat{P}^{a,b}_{R,j}(s)+ \widehat{C}_j^{a,b}\widehat{P}^{a,b}_{R,j+1}(s)\right)\\
&+\frac{A_j^{a,b}}{\Gamma(\alpha)}(1-x)^{\alpha}\left[\widehat{A}_j^{a,b}P^{a,b}_{j-1}(1)
+\widehat{B}_j^{a,b}P^{a,b}_{j}(1)+ \widehat{C}_j^{a,b}P^{a,b}_{j+1}(1)\right].
 \end{aligned}
\end{equation}
Note that $\widehat{P}_{R,0}^{a,b,\alpha}(x)$ and $\widehat{P}_{R,1}^{a,b,\alpha}(x)$ can
be easily calculated.  Hence,  we derive the recurrence relation from \eqref{eq:fint1-3}
below
\begin{equation} \label{eq:fint2}
\left\{\begin{aligned}
&\widehat{P}_{R,0}^{a,b,\alpha}(x)=\frac{(1-x)^{\alpha}}{\Gamma(\alpha+1)},\\
&\widehat{P}_{R,1}^{a,b,\alpha}(x)=\frac{a+b+2}{2}\left(\frac{x(1-x)^{\alpha}}{\Gamma(\alpha+1)}
+\frac{\alpha(1-x)^{\alpha+1}}{\Gamma(\alpha+2)}\right)+\frac{a-b}{2}\widehat{P}_{R,0}^{a,b,\alpha}(x),\\
&\widehat{P}_{R,j+1}^{a,b,\alpha}(x)=\frac{A_j^{a,b}x-(B_j^{a,b}+{\alpha}A_j^{a,b}\widehat{B}_j^{a,b})}
{1+{\alpha}A_j^{a,b}\widehat{C}_j^{a,b}}\widehat{P}_{R,j}^{a,b,\alpha}(x)
-\frac{C_j^{a,b}+{\alpha}A_j^{a,b}\widehat{A}_j^{a,b}}{1+{\alpha}A_j^{a,b}\widehat{C}_j^{a,b}}
\widehat{P}_{R,j-1}^{a,b,\alpha}(x)\\
&{\quad}+\frac{{\alpha}A_j^{a,b}\left(\widehat{A}_j^{a,b}P^{a,b}_{j-1}(1)+\widehat{B}_j^{a,b}P^{a,b}_{j}(1)
+\widehat{C}_j^{a,b}P^{a,b}_{j+1}(1)\right)}
{\Gamma(\alpha+1)\left(1+{\alpha}A_j^{a,b}\widehat{C}_j^{a,b}\right)}
(1-x)^{\alpha},{\quad\quad}j{\,\geq\,}1.
\end{aligned}\right.
\end{equation}

If $a=b=0$, the recurrence formulas  \eqref{eq:fint1} and \eqref{eq:fint2} are reduced to
\begin{equation} \label{eq:fint3}
\left\{\begin{aligned}
\widehat{L}_{L,0}^{\alpha}(x)&=\frac{(x+1)^{\alpha}}{\Gamma(\alpha+1)},\quad
\widehat{L}_{L,1}^{\alpha}(x)=\frac{x(x+1)^{\alpha}}{\Gamma(\alpha+1)}
-\frac{\alpha(x+1)^{\alpha+1}}{\Gamma(\alpha+2)},\\
\widehat{L}_{L,j+1}^{\alpha}(x)&=\frac{1}{j+1+\alpha}\Big\{(2j+1)x\widehat{L}_{L,j}^{\alpha}(x)
-(j-\alpha)\widehat{L}_{L,j-1}^{\alpha}(x)\Big\},{\quad}j{\,\geq\,}1
\end{aligned}\right.
\end{equation}
and
\begin{equation} \label{eq:fint4}
\left\{\begin{aligned}
\widehat{L}_{R,0}^{\alpha}(x)&=\frac{(1-x)^{\alpha}}{\Gamma(\alpha+1)},\quad
\widehat{L}_{R,1}^{\alpha}(x)=\frac{x(1-x)^{\alpha}}{\Gamma(\alpha+1)}
+\frac{\alpha(1-x)^{\alpha+1}}{\Gamma(\alpha+2)},\\
\widehat{L}_{R,j+1}^{\alpha}(x)&=\frac{1}{j+1+\alpha}\Big\{(2j+1)x\widehat{L}_{R,j}^{\alpha}(x)
-(j-\alpha)\widehat{L}_{R,j-1}^{\alpha}(x)\Big\},{\quad}j{\,\geq\,}1,
\end{aligned}\right.
\end{equation}
respectively, where $\widehat{L}_{L,j}^{\alpha}(x)=\widehat{P}_{L,j}^{0,0,\alpha}(x)$ and
$\widehat{L}_{R,j}^{\alpha}(x)=\widehat{P}_{R,j}^{0,0,\alpha}(x)$.

For $a=b=-1/2$, from  \eqref{eq:fint1}, \eqref{eq:fint2},  and  the relation
$P^{-\frac{1}{2},-\frac{1}{2}}_{j}(x)=\frac{\Gamma(j+1/2)}{j!\sqrt{\pi}}T_j(x)$, one can obtain
\begin{equation} \label{eq:fint5}
\left\{\begin{aligned}
\widehat{T}_{L,0}^{\alpha}(x)=&\frac{(x+1)^{\alpha}}{\Gamma(\alpha+1)},{\qquad}
\widehat{T}_{L,1}^{\alpha}(x)=\frac{x(x+1)^{\alpha}}{\Gamma(\alpha+1)}
-\frac{\alpha(x+1)^{\alpha+1}}{\Gamma(\alpha+2)},\\
\widehat{T}_{L,2}^{\alpha}(x)=&\frac{2(1+x)^{2+\alpha}}{\Gamma(3+\alpha)}
-\frac{4(1+x)^{1+\alpha}}{\Gamma(2+\alpha)}+\frac{(1+x)^{\alpha}}{\Gamma(1+\alpha)},\\
\widehat{T}_{L,j+1}^{\alpha}(x)=&\frac{2(j+1)x}{j+1+\alpha}\widehat{T}_{L,j}^{\alpha}(x)
-\frac{(j+1)(j-1-\alpha)}{(j+1+\alpha)(j-1)}\widehat{T}_{L,j-1}^{\alpha}(x)\\
&+\frac{2(-1)^{j}\alpha(x+1)^{\alpha}}{\Gamma(\alpha+1)(j+1+\alpha)(j-1)},\quad
j{\,\geq\,}2
\end{aligned}\right.
\end{equation}
and
\begin{equation} \label{eq:fint6}
\left\{\begin{aligned}
\widehat{T}_{R,0}^{\alpha}(x)=&\frac{(1-x)^{\alpha}}{\Gamma(\alpha+1)},{\qquad}
\widehat{T}_{R,1}^{\alpha}(x)=\frac{x(1-x)^{\alpha}}{\Gamma(\alpha+1)}
+\frac{\alpha(1-x)^{\alpha+1}}{\Gamma(\alpha+2)},\\
\widehat{T}_{R,2}^{\alpha}(x)=&\frac{2(1-x)^{2+\alpha}}{\Gamma(3+\alpha)}
-\frac{4(1-x)^{1+\alpha}}{\Gamma(2+\alpha)}+\frac{(1-x)^{\alpha}}{\Gamma(1+\alpha)},\\
\widehat{T}_{R,j+1}^{\alpha}(x)=&\frac{2(j+1)x}{j+1+\alpha}\widehat{T}_{R,j}^{\alpha}(x)
-\frac{(j+1)(j-1-\alpha)}{(j+1+\alpha)(j-1)}\widehat{T}_{R,j-1}^{\alpha}(x)\\
&+\frac{2\alpha(1-x)^{\alpha}}{\Gamma(\alpha+1)(j+1+\alpha)(j-1)},\quad
j{\,\geq\,}2,
\end{aligned}\right.
\end{equation}
respectively,
where $\widehat{T}_{L,j}^{\alpha}(x)=\frac{1}{\Gamma(\alpha)}\int_{-1}^x(x-s)^{\alpha-1}T_{j}(s)\dx[s]$
and $\widehat{T}_{R,j}^{\alpha}(x)=\frac{1}{\Gamma(\alpha)}\int_{x}^1(s-x)^{\alpha-1}T_{j}(s)\dx[s]$.
One can refer to \cite{LiZengLiu2012},
in which the detailed derivations of \eqref{eq:fint3} and \eqref{eq:fint5} can be found.

Let $u(x)$ be a function defined on the interval $[-1,1]$ and $N$
be a positive integer.   Denote $x_j(j=0,1,...,N)$ as the Jacobi--Gauss--Lobatto (JGL)
points defined on the interval $[-1,1]$. Then the JGL interpolation
of $u(x)$ is given by
\begin{equation}\label{eq:in1}
I_Nu(x)=\sum_{j=0}^Nu(x_j)F_j(x)=\sum_{j=0}^N\tilde{p}_jP_j^{a,b}(x),
\end{equation}
where $F_j(x)$ is the Lagrange base function.

In \cite{LiZengLiu2012}, authors obtained the formulas for numerically calculating
$D_{-1,x}^{-\alpha}I_Nu(x)$ and  ${}_CD^{\alpha}_{-1,x}I_Nu(x)$ on the JGL
points $\{x_j\}$.  $D_{-1,x}^{-\alpha}I_Nu(x)(\alpha>0)$ at $x=x_j,j=0,1,...,N$  can be calculated by
the following formula \cite{LiZengLiu2012}
\begin{equation}\label{eq:fint7}
\left(\begin{array}{c}
  D_{-1,x_0}^{-\alpha}I_Nu(x_0) \\
  D_{-1,x_1}^{-\alpha}I_Nu(x_1) \\
    \vdots\\
 D_{-1,x_N}^{-\alpha}I_Nu(x_N)
\end{array}\right)={\widehat{D}}^{(a,b,\alpha)}_L(\tilde{p}_0,\tilde{p}_1,\cdots,\tilde{p}_N)^T,
\end{equation}
where $D_{-1,x_i}^{-\alpha}I_Nu(x_i)=\left[D_{-1,x}^{-\alpha}I_Nu\right]_{x=x_i}$ and
${\widehat{D}}^{(a,b,\alpha)}_L \in \mathbb{R}^{(N+1)\times(N+1)}$ satisfying
$$\left({\widehat{D}}^{(a,b,\alpha)}_L\right)_{i,j}=\widehat{P}_{L,j}^{a,b,\alpha}({x}_i),{\quad}i,j=0,...,N.$$
The computation of ${}_C{D}^{\alpha}_{-1,x}I_Nu(x)(n-1<\alpha<n,n\in Z^+)$ at $x=x_j,j=0,1,...,N$ is given by
\cite{LiZengLiu2012}
\begin{equation}\label{eq:fint8}
\left(\begin{array}{c}
  D_{-1,x_0}^{\alpha}I_Nu(x_0) \\
  D_{-1,x_1}^{\alpha}I_Nu(x_1) \\
    \vdots\\
 D_{-1,x_N}^{\alpha}I_Nu(x_N)
\end{array}\right)={}_C{\widehat{D}}_L^{(a,b,\alpha)}(\tilde{p}_0,\tilde{p}_1,\cdots,\tilde{p}_N)^T,
\end{equation}
where ${}_{C}D_{-1,x_i}^{\alpha}I_Nu(x_i)=\left[{}_{C}D_{-1,x}^{\alpha}I_Nu\right]_{x=x_i}$ and
${}_C{\widehat{D}}^{(a,b,\alpha)}_L \in \mathbb{R}^{(N+1)\times(N+1)}$ satisfying
$$\left({}_C{\widehat{D}}_L^{(a,b,\alpha)}\right)_{i,j}=d_{j,n}^{a,b}\,\widehat{P}^{a+n,b+n,n-\alpha}_{L,j}(x_i),
{\quad}i,j=0,...,N$$
with $d_{j,n}^{a,b}$ given by \eqref{jacobi3.2}.

Suppose that the Lagrange base function $F_j(x)$  can be expressed as
\begin{equation}\label{eq:fint9}
F_j(x)=\sum_{k=0}^Nc_{k,j}P_k^{a,b}(x).
\end{equation}
Then $c_{k,j}$ can be determined by the following relation \cite{ShenTangWang2011}
\begin{equation} \label{eq:ckj}
c_{k,j}=\left\{\begin{aligned}
&\frac{P_k^{a,b}(x_j)\omega_j}{\gamma_k^{a,b}},{\,\,\,\,\qquad}k=0,1,...,N-1,\\
&\frac{P_N^{a,b}(x_j)\omega_j}{(2+\frac{a+b+1}{N})\gamma_N^{a,b}},{\quad}k=N,
\end{aligned}\right.
\end{equation}
in which $\gamma_k^{a,b}$ is defined by \eqref{jacobi1-3}, $x_j$ is the JGL point on $[-1,1]$.

From \eqref{eq:in1},  \eqref{eq:fint9}, and \eqref{eq:ckj}, we obtain
\begin{equation} \label{eq:fint10}
\begin{aligned}
\left(
\begin{array}{c}
 \tilde{p}_0^{a,b} \\
 \tilde{p}_1^{a,b} \\
   \vdots \\
 \tilde{p}_N^{a,b} \\
    \end{array}
     \right)=\left(
\begin{array}{cccc}
    c_{0,0} &c_{0,1} &\cdots &c_{0,N}  \\
    c_{1,0} &c_{1,1} &\cdots &c_{1,N}  \\
   \vdots   &\vdots &\ddots &\vdots  \\
    c_{N,0} &c_{N,1} &\cdots &c_{N,N}  \\
    \end{array}
    \right)\left(
\begin{array}{c}
 u(x_0) \\
 u(x_1) \\
   \vdots \\
 u(x_N) \\
    \end{array}
     \right)=M(u(x_0),u(x_1),\cdots,u(x_N))^T.
\end{aligned}
\end{equation}

Hence, the fractional integral $D_{-1,x}^{-\alpha}I_Nu(x)$ at $x=x_j$ defined by \eqref{eq:fint7}
can be written into the following equivalent form
\begin{equation}\label{eq:fint11}
\left(\begin{array}{c}
  D_{-1,x_0}^{-\alpha}I_Nu(x_0) \\
  D_{-1,x_1}^{-\alpha}I_Nu(x_1) \\
    \vdots\\
 D_{-1,x_N}^{-\alpha}I_Nu(x_N)
\end{array}\right)=\left({\widehat{D}}^{(a,b,\alpha)}_LM\right)(u(x_0),u(x_1),\cdots,u(x_N))^T,
\end{equation}
where the matrices ${D}^{(a,b,\alpha)}_L$ and $M$ are defined as in \eqref{eq:fint7}
and \eqref{eq:fint10}, respectively.
We can similarly  obtain the equivalent form of \eqref{eq:fint8} as follows
\begin{equation}\label{eq:fint12}
\left(\begin{array}{c}
  {}_{C}D_{-1,x_0}^{\alpha}I_Nu(x_0) \\
  {}_{C}D_{-1,x_1}^{\alpha}I_Nu(x_1) \\
    \vdots\\
 {}_{C}D_{-1,x_N}^{\alpha}I_Nu(x_N)
\end{array}\right)=\left({}_C{\widehat{D}}^{(a,b,\alpha)}_LM\right)(u(x_0),u(x_1),\cdots,u(x_N))^T,
\end{equation}
where the matrices ${}_C{D}^{(a,b,\alpha)}_L$ and $M$ are defined as in \eqref{eq:fint8}
and \eqref{eq:fint10}, respectively.

Next, we define the   matrix ${}_C{D}^{(a,b,\alpha)}_L,\alpha \in \mathbb{R}$ as follows:
\begin{equation} \label{eq:fracmat}
{}_C{D}^{(a,b,\alpha)}_L=\left\{\begin{aligned}
&{\widehat{D}}^{(a,b,-\alpha)}_LM,{\quad}{\quad}  M \text{ and } {\widehat{D}}^{(a,b,-\alpha)}_L
\text{ are defined as in \eqref{eq:fint11} for } \alpha<0, \\
&{}_C{\widehat{D}}^{(a,b,\alpha)}_LM,{\quad} {\quad}  M \text{ and } {}_C{\widehat{D}}^{(a,b,\alpha)}_L
\text{ are defined as in \eqref{eq:fint12} for }\alpha>0.
\end{aligned}\right.
\end{equation}

We call the matrix ${}_{C}D_L^{(a,b,\alpha)}(\alpha>0)$ the fractional differential matrix with respect to the
left Caputo derivative operator. We can similarly obtain the corresponding fractional differential matrix
 ${}_{C}D_R^{(a,b,\alpha)}$ with respect to the right fractional integral and Caputo derivative operators.

Since the left Caputo and Riemann--Liouville derivative operators have the following relation
\begin{equation} \label{eq:rl-caputo}
{}_{RL}D_{x_0,x}^{\alpha}f(x)={}_{C}D_{x_0,x}^{\alpha}f(x)
  +\sum_{k=0}^{n-1}\frac{f^{(k)}(x_0)(x-x_0)^{k-\alpha}}{\Gamma(k+1-\alpha)},
  {\quad}n-1<\alpha<n,n \in \mathbb{N}.
\end{equation}
Hence, we can obtain
\begin{equation}\label{eq:fint13}
\left(\begin{array}{c}
  {}_{RL}D_{-1,x_0}^{\alpha}I_Nu(x_0) \\
  {}_{RL}D_{-1,x_1}^{\alpha}I_Nu(x_1) \\
    \vdots\\
 {}_{RL}D_{-1,x_N}^{\alpha}I_Nu(x_N)
\end{array}\right)=\left({}_{RL}{{D}}^{(a,b,\alpha)}_L\right)(u(x_0),u(x_1),\cdots,u(x_N))^T,
\end{equation}
where  ${}_{RL}D_{-1,x_i}^{\alpha}I_Nu(x_i)=\left[{}_{RL}D_{-1,x}^{\alpha}I_Nu\right]_{x=x_i}$ and
the matrix ${}_{RL}{{D}}^{(a,b,\alpha)}_L  \in \mathbb{R}^{(N+1)\times(N+1)}$ satisfying
\begin{equation}\label{eq:fint14}
\left({}_{RL}{{D}}^{(a,b,\alpha)}_L\right)_{i,j}
=\left({}_{C}{{D}}^{(a,b,\alpha)}_L\right)_{i,j}
+\sum_{k=0}^{n-1}\frac{D^{(k)}_{0,j}(x_{i}+1)^{k-\alpha}}{\Gamma(k+1-\alpha)},{\quad}
i,j=0,1,\cdots,N,
\end{equation}
in which ${}_{RL}D_{-1,x_i}^{\alpha}I_Nu(x_i)=\left[{}_{RL}D_{-1,x}^{\alpha}I_Nu\right]_{x=x_i}$ and
$D^{(k)}={}_C{D}^{(a,b,k)}_L$ is the $k$th-order classical differential matrix,
see \cite{ShenTangWang2011} for more information.
We can similarly obtain
\begin{equation}\label{eq:fint15}
\left(\begin{array}{c}
  {}_{RL}D_{x_0,1}^{\alpha}I_Nu(x_0) \\
  {}_{RL}D_{x_1,1}^{\alpha}I_Nu(x_1) \\
    \vdots\\
 {}_{RL}D_{x_{N},1}^{\alpha}I_Nu(x_{N})
\end{array}\right)=\left({}_{RL}{{D}}^{(a,b,\alpha)}_R\right)(u(x_0),u(x_1),\cdots,u(x_N))^T,
\end{equation}
where ${}_{RL}D_{x_i,1}^{\alpha}I_Nu(x_i)=\left[{}_{RL}D_{x,1}^{\alpha}I_Nu\right]_{x=x_i}$ and
the matrix ${}_{RL}{{D}}^{(a,b,\alpha)}_R  \in \mathbb{R}^{(N+1)\times(N+1)}$ satisfying
\begin{equation}\label{eq:fint16}
\left({}_{RL}{{D}}^{(a,b,\alpha)}_R\right)_{i,j}
=\left({}_{C}{{D}}^{(a,b,\alpha)}_R\right)_{i,j}
+\sum_{k=0}^{n-1}\frac{(-1)^kD^{(k)}_{N,j}(1-x_{i})^{k-\alpha}}{\Gamma(k+1-\alpha)},
{\quad}i,j=0,1,\cdots,N.
\end{equation}
\begin{remark}
Note that the first row of the matrix ${}_{RL}{{D}}^{(a,b,\alpha)}_L$ and
the last row of the matrix ${}_{RL}{{D}}^{(a,b,\alpha)}_R$ have no sense, except
that $\alpha$ is a positive integer, i.e., $\alpha=n,n\in \mathbb{N}$, and we set
${}_{RL}{{D}}^{(a,b,n)}_L={}_{C}{{D}}^{(a,b,n)}_L$ and
${}_{RL}{{D}}^{(a,b,n)}_R={}_{C}{{D}}^{(a,b,n)}_R$.
\end{remark}

Next, we investigate the eigenvalues of the fractional differential operators. Consider first the
following model problem
\begin{equation}\label{eq:fint17}
D^{\alpha}u(x)=\lambda u(x),x\in (-1,1),{\quad}u(-1)=u(1)=0, {\quad}1<\alpha{\,\leq\,}2.
\end{equation}
where $D^{\alpha}$ denotes the left (or right) Caputo derivative operator ${}_{C}D_{-1,x}^{\alpha}$
(or  ${}_{C}D_{x,1}^{\alpha}$), the left (or right) Riemann--Liouville derivative operator
${}_{RL}D_{-1,x}^{\alpha}$ (or  ${}_{RL}D_{x,1}^{\alpha}$),  the Riesz-Caputo fractional derivative
operator ${}^C_{RZ}D_x^{\alpha}$, or the Riesz fractional derivative
operator ${}_{RZ}D_x^{\alpha}$.

For simplicity, we first consider the case  $D^{\alpha}={}_{C}D_{-1,x}^{\alpha}$ in \eqref{eq:fint17}.
Denote  ${x}_i(i=0,1,...,N)$  as the JGL points on the interval $[-1,1]$.
Suppose that $u(x)$ can be approximated by
\begin{equation}\label{eq:fint17-2}
u(x) \approx I_Nu(x)=\sum_{j=0}^Nu_jF_j(x).
\end{equation}
Replacing $u(x)$ with $I_Nu(x)$ in \eqref{eq:fint17} and letting $x={x}_i$ yield
\begin{equation}\label{eq:fint18}
\left[{}_{C}D_{-1,x}^{\alpha}I_Nu(x)\right]_{x={x}_i}=\lambda (I_Nu)({x}_i),{\quad}i=1,2,...,N-1.
\end{equation}
Using the boundary conditions $(I_Nu)(-1)=(I_Nu)(1)=0$, we obtain  the matrix representation
of \eqref{eq:fint18} as follows
\begin{equation}\label{eq:fint19}
M_{C,L}^{(\alpha,a,b)}\mathbf{u}=\lambda \mathbf{u},
\end{equation}
where $\mathbf{u}=(u_1,u_2,...,u_{N-1})^T$ and $M_{C,L}^{(\alpha,a,b)} \in \mathbb{R}^{(N-1)\times(N-1)}$
satisfying
$$\left(M_{C,L}^{(\alpha,a,b)} \right)_{i,j}=\left({}_{C}{{D}}^{(a,b,\alpha)}_L\right)_{i+1,j+1},
{\quad}i,j=0,1,...,N-2.$$

If $D^{\alpha}$  in \eqref{eq:fint17} is chosen as  $D^{\alpha}={}_{C}D_{x,1}^{\alpha}$,
$D^{\alpha}={}_{RL}D_{-1,x}^{\alpha}$,   $D^{\alpha}={}_{RL}D_{x,1}^{\alpha}$,
$D^{\alpha}={}^C_{RZ}D_{x}^{\alpha}$, or
$D^{\alpha}={}_{RZ}D_{x}^{\alpha}$, then we can similarly derive the linear system as \eqref{eq:fint19}
with the coefficient matrices denoted by  $M_{C,R}^{(\alpha,a,b)}$, $M_{RL,L}^{(\alpha,a,b)}$,
$M_{RL,R}^{(\alpha,a,b)}$, $M_{C,RZ}^{(\alpha,a,b)}$, or $M_{RZ}^{(\alpha,a,b)}$.
Denote by
$M_1^{(\alpha,a,b)}=M_{C,L}^{(\alpha,a,b)}$, $M_2^{(\alpha,a,b)}=M_{C,R}^{(\alpha,a,b)}$,
$M_3^{(\alpha,a,b)}=M_{RL,L}^{(\alpha,a,b)}$,
$M_4^{(\alpha,a,b)}=M_{RL,R}^{(\alpha,a,b)}$,   $M_5^{(\alpha,a,b)}=M_{RC}^{(\alpha,a,b)}$,
and $M_6^{(\alpha,a,b)}=M_{RZ}^{(\alpha,a,b)}$.
It is well known that  the spectral radius  $\rho\left(M_{i}^{(2,0,0)}\right)$
of  $M_{i}^{(2,0,0)}(i=1,2,3,4,5,6)$  satisfy the following
relation \cite{ShenTangWang2011}
$$\rho\left(M_{i}^{(2,0,0)}\right){\,\leq\,}C_0N^{4}, {\quad}C_0>0 \text{ is independent of } N.$$

Is it possible that the spectral radius  $\rho\left(M_{i}^{(\alpha,a,b)}\right)$ of
$M_{i}^{(\alpha,a,b)}$ satisfies the following relation?
\begin{equation}\label{eq:fint20}
\rho\left(M_{i}^{(\alpha,a,b)}\right){\,\leq\,}C_0N^{2\alpha}, {\quad}C_0>0 \text{ is independent of } N.
\end{equation}

\begin{figure}
\begin{center}
\begin{minipage}{0.4\textwidth}\centering
\includegraphics[scale=0.35]{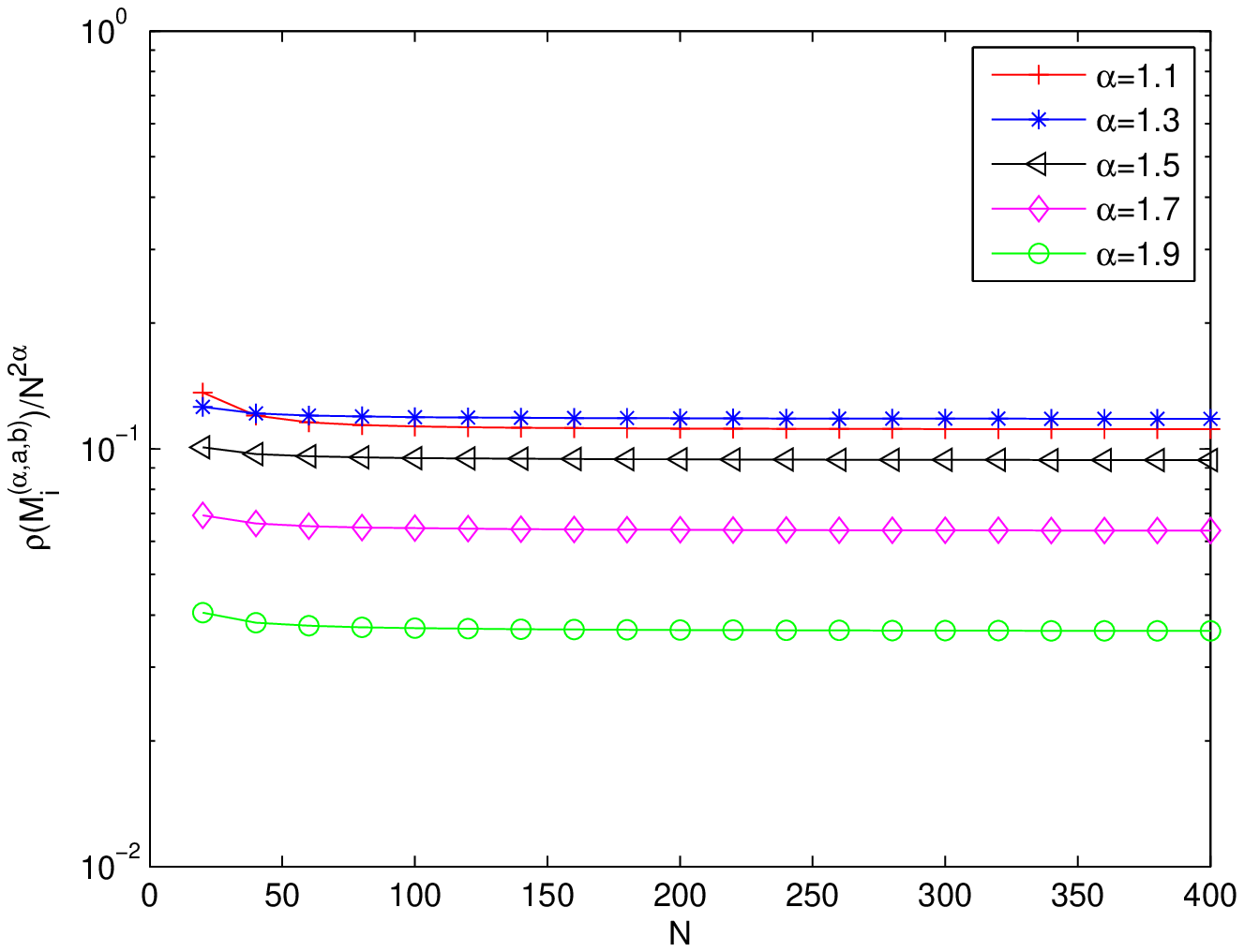} \par {(a)  $M_{i}^{(\alpha,a,b)}=M_{C,L}^{(\alpha,0,0)}$.}
\end{minipage}
\begin{minipage}{0.4\textwidth}\centering
\includegraphics[scale=0.35]{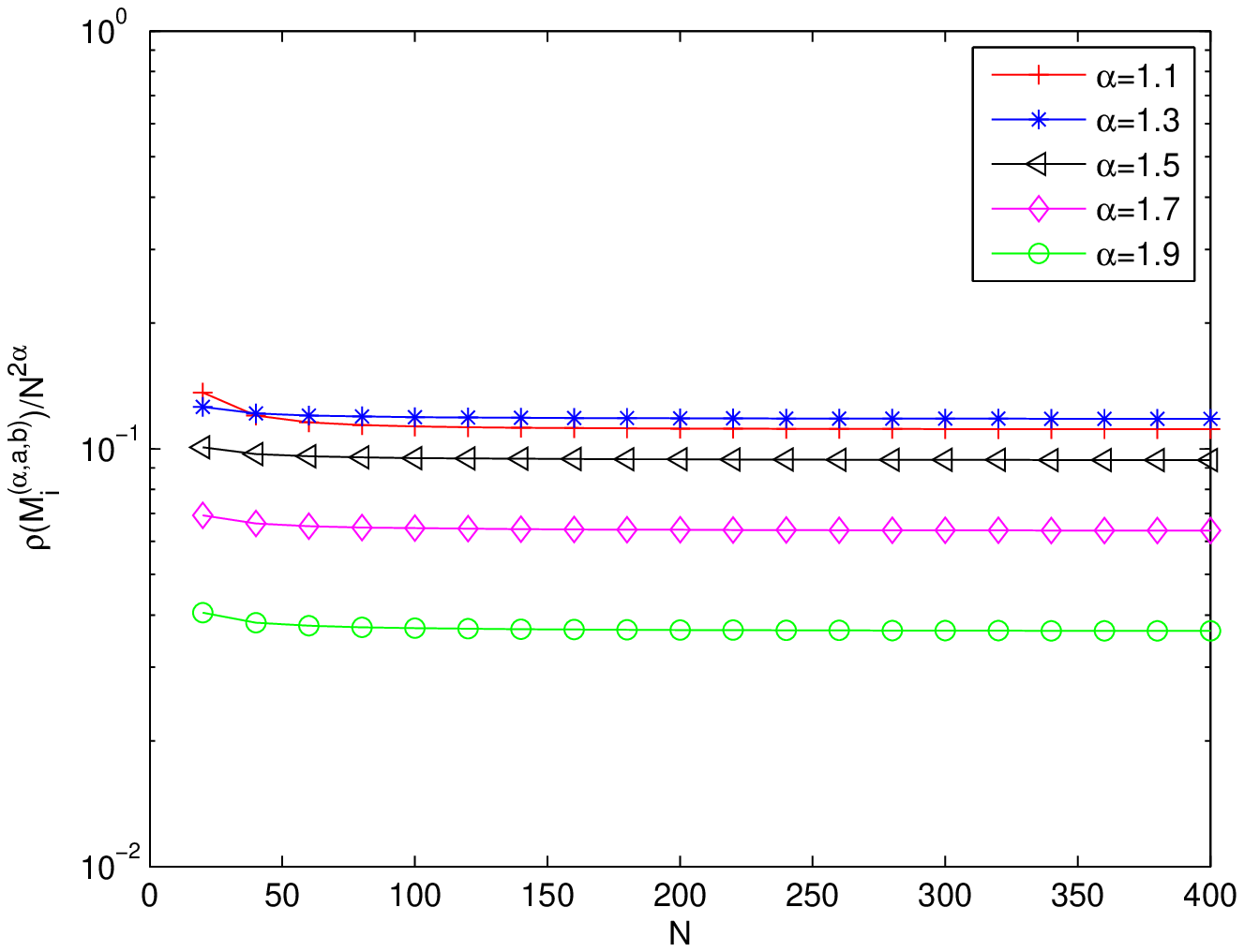} \par {(b) $M_{i}^{(\alpha,a,b)}=M_{C,R}^{(\alpha,0,0)}$.}
\end{minipage}\\
\begin{minipage}{0.4\textwidth}\centering
\includegraphics[scale=0.35]{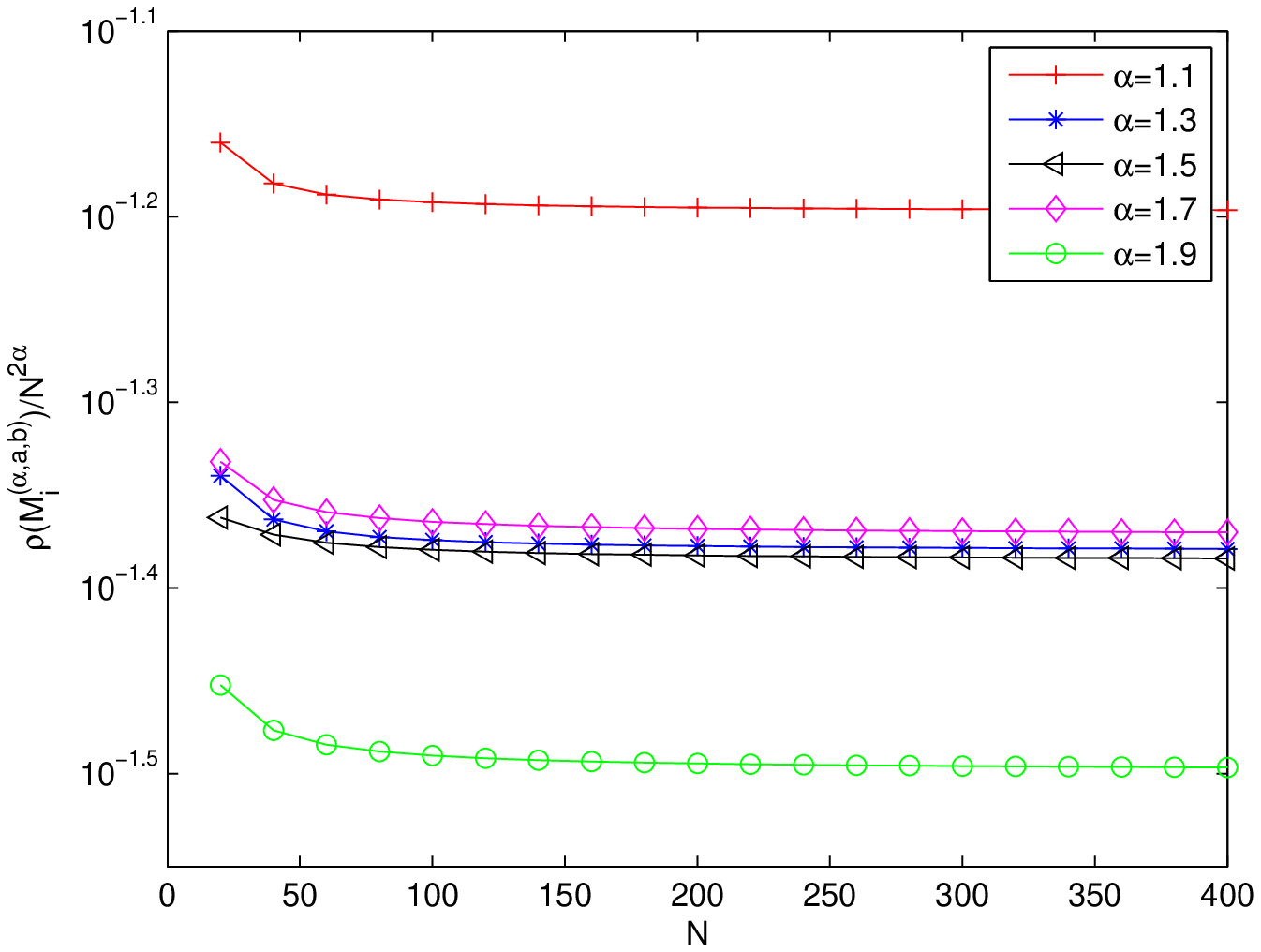} \par {(c)  $M_{i}^{(\alpha,a,b)}=M_{RL,L}^{(\alpha,0,0)}$.}
\end{minipage}
\begin{minipage}{0.4\textwidth}\centering
\includegraphics[scale=0.35]{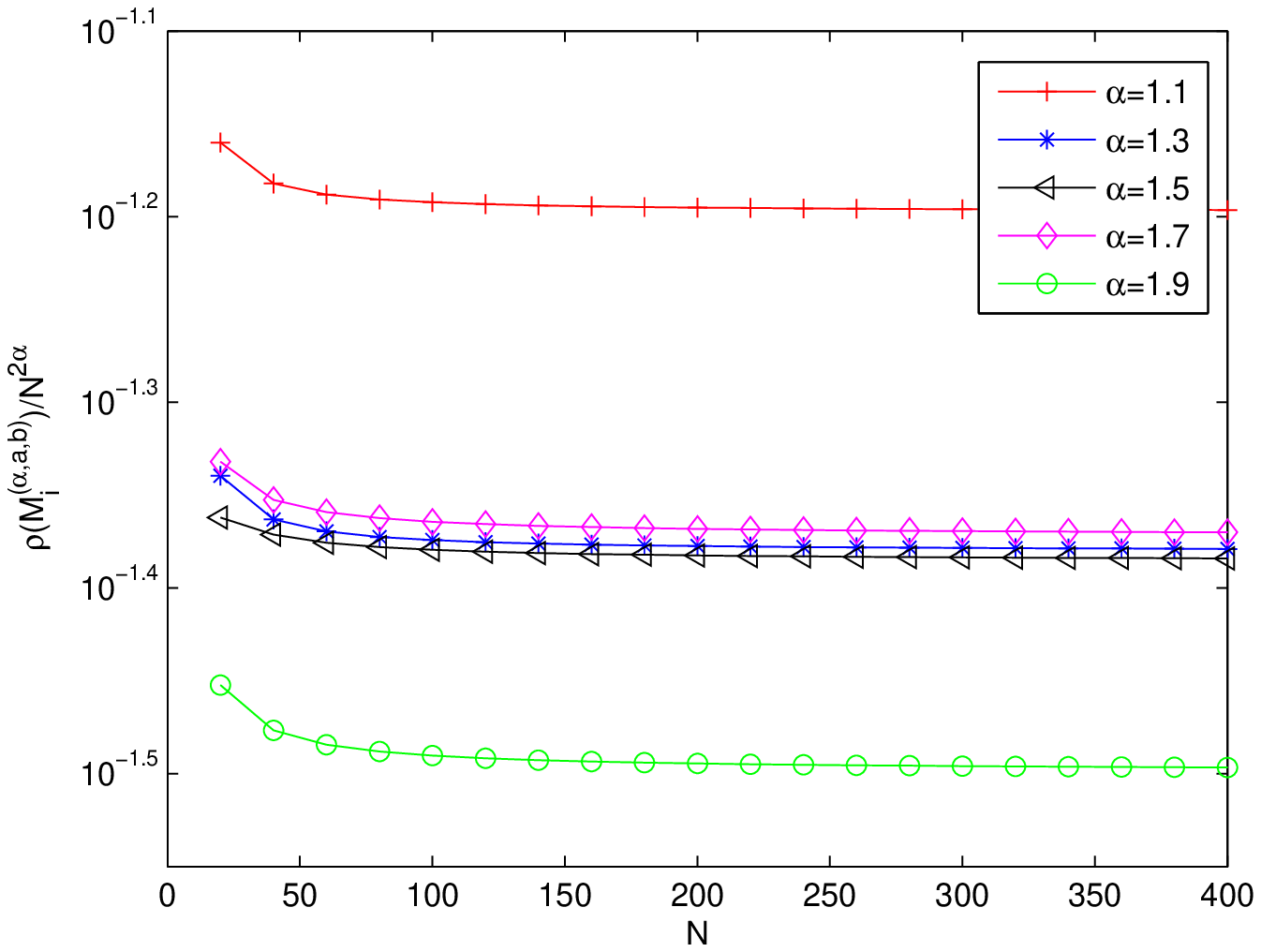} \par {(d)  $M_{i}^{(\alpha,a,b)}=M_{RL,R}^{(\alpha,0,0)}$.}
\end{minipage}
\begin{minipage}{0.4\textwidth}\centering
\includegraphics[scale=0.35]{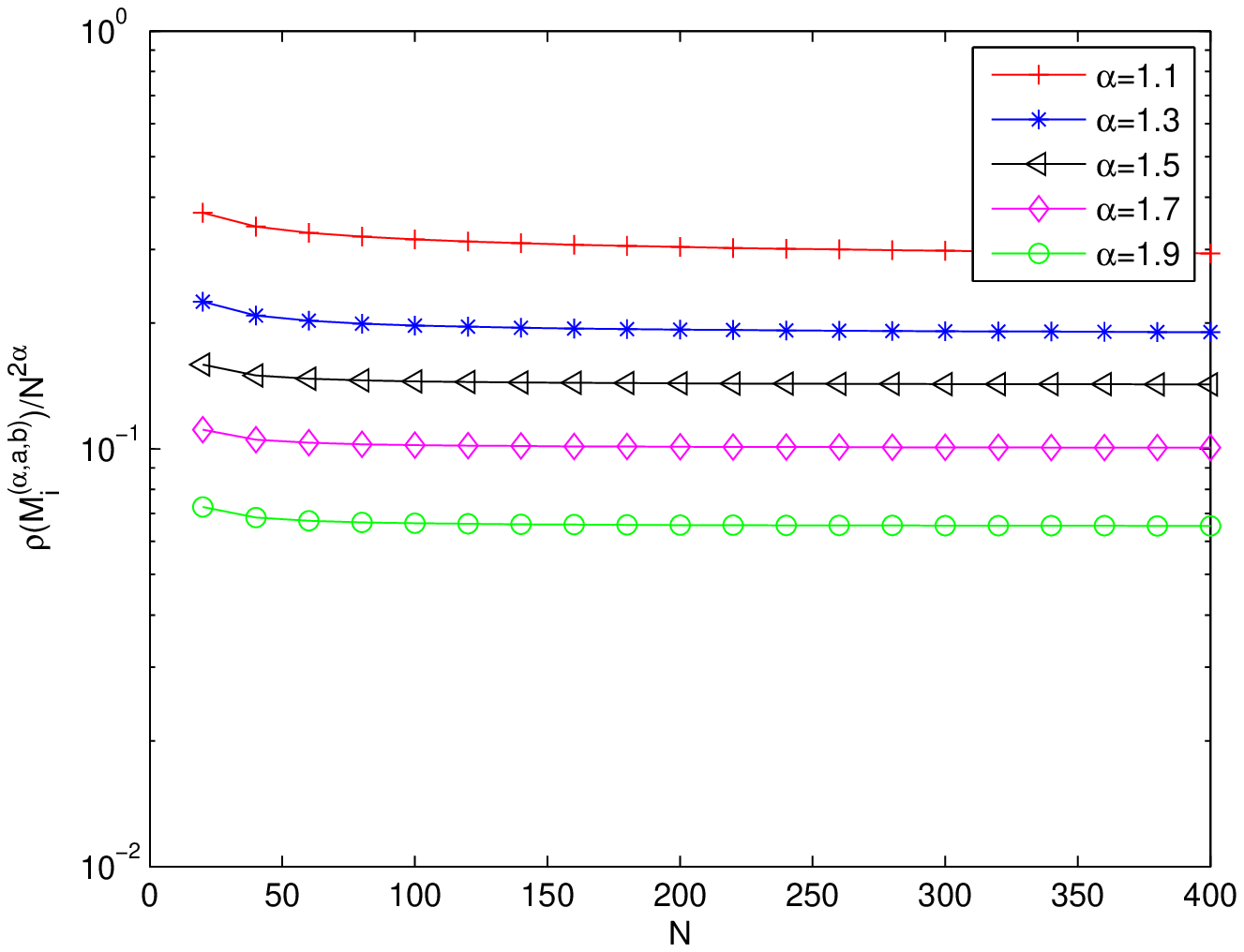} \par {(e) $M_{i}^{(\alpha,a,b)}=M_{RC}^{(\alpha,0,0)}$.  }
\end{minipage}
\begin{minipage}{0.4\textwidth}\centering
\includegraphics[scale=0.35]{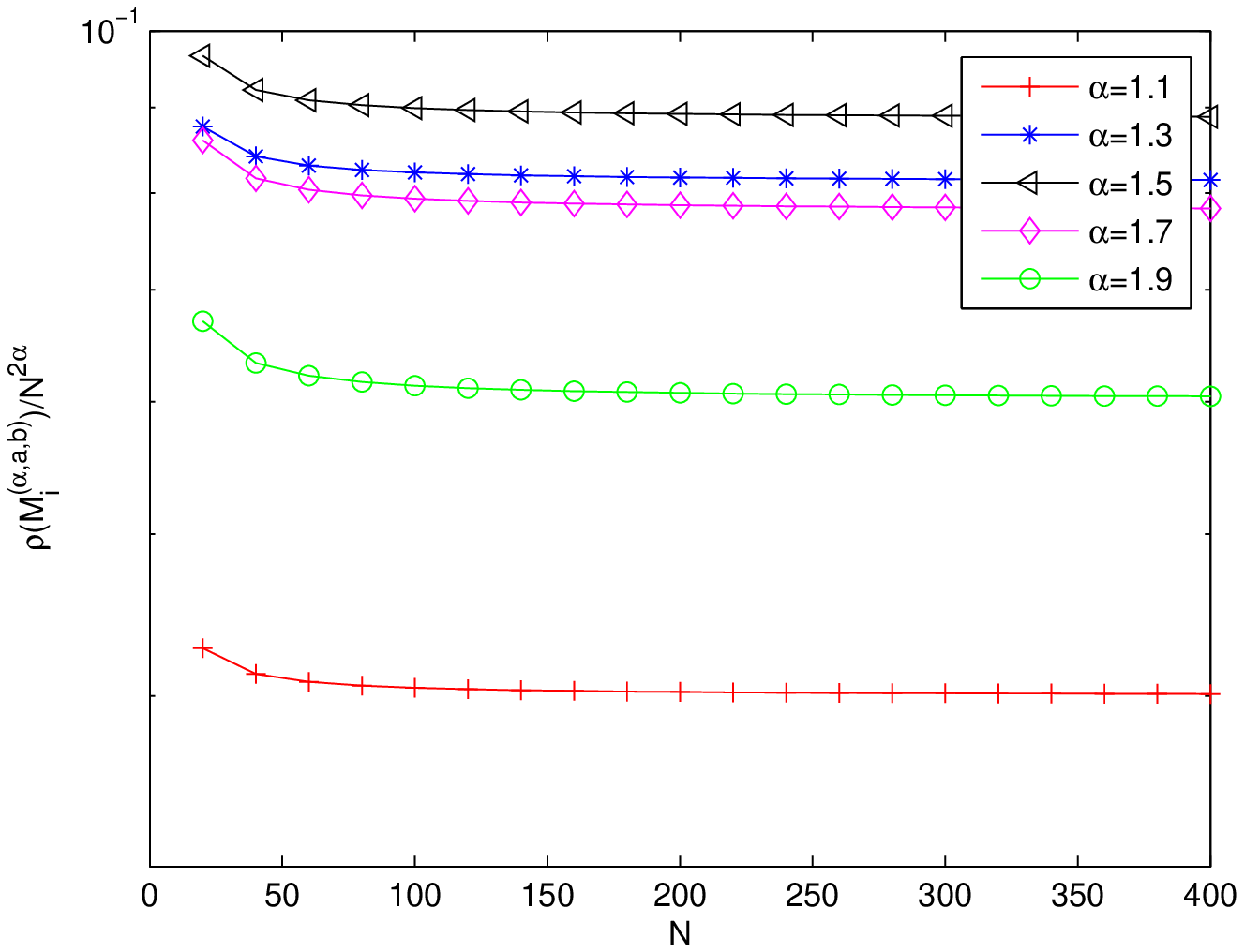} \par {(f) $M_{i}^{(\alpha,a,b)}=M_{RZ}^{(\alpha,0,0)}$. }
\end{minipage}
\end{center}
\caption{The boundedness of $\rho\left(M_{i}^{(\alpha,a,b)}\right)/N^{2\alpha}$ for $a=b=0$.\label{fig1}}
\end{figure}

\begin{figure}
\begin{center}
\begin{minipage}{0.4\textwidth}\centering
\includegraphics[scale=0.35]{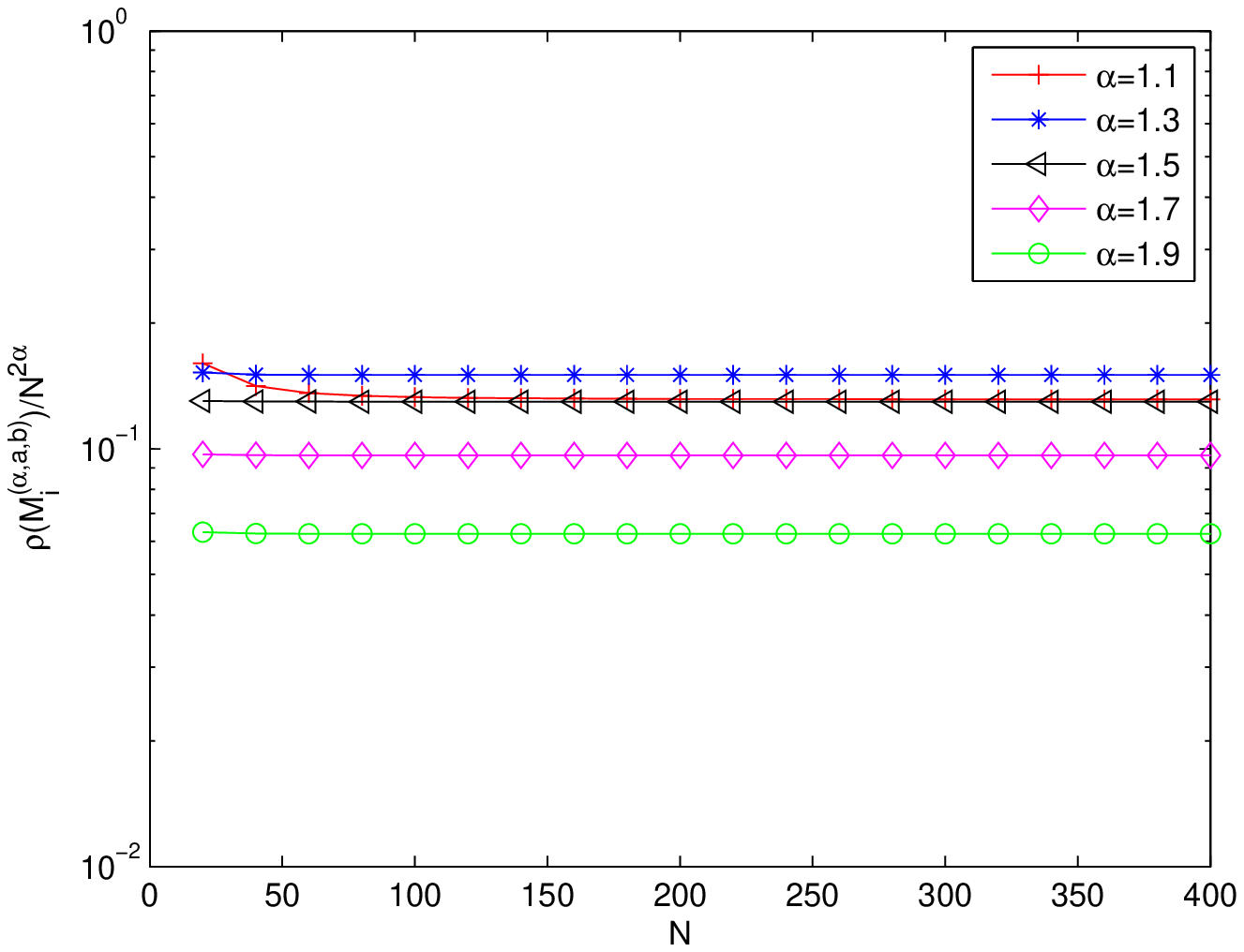}
\par {(a)  $M_{i}^{(\alpha,a,b)}=M_{C,L}^{(\alpha,-\frac{1}{2},-\frac{1}{2})}$.}
\end{minipage}
\begin{minipage}{0.4\textwidth}\centering
\includegraphics[scale=0.35]{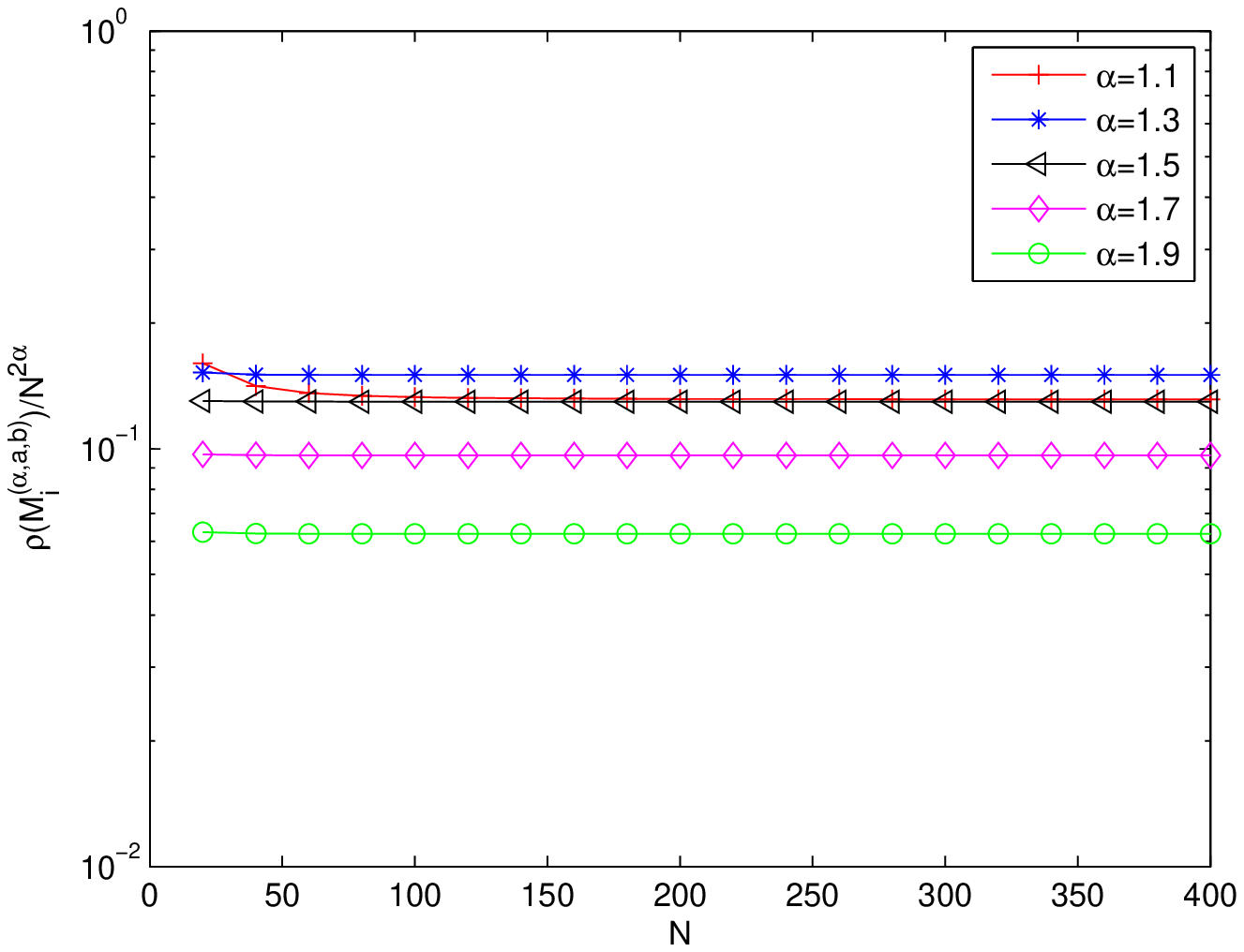}
\par {(b) $M_{i}^{(\alpha,a,b)}=M_{C,R}^{(\alpha,-\frac{1}{2},-\frac{1}{2})}$.}
\end{minipage}\\
\begin{minipage}{0.4\textwidth}\centering
\includegraphics[scale=0.35]{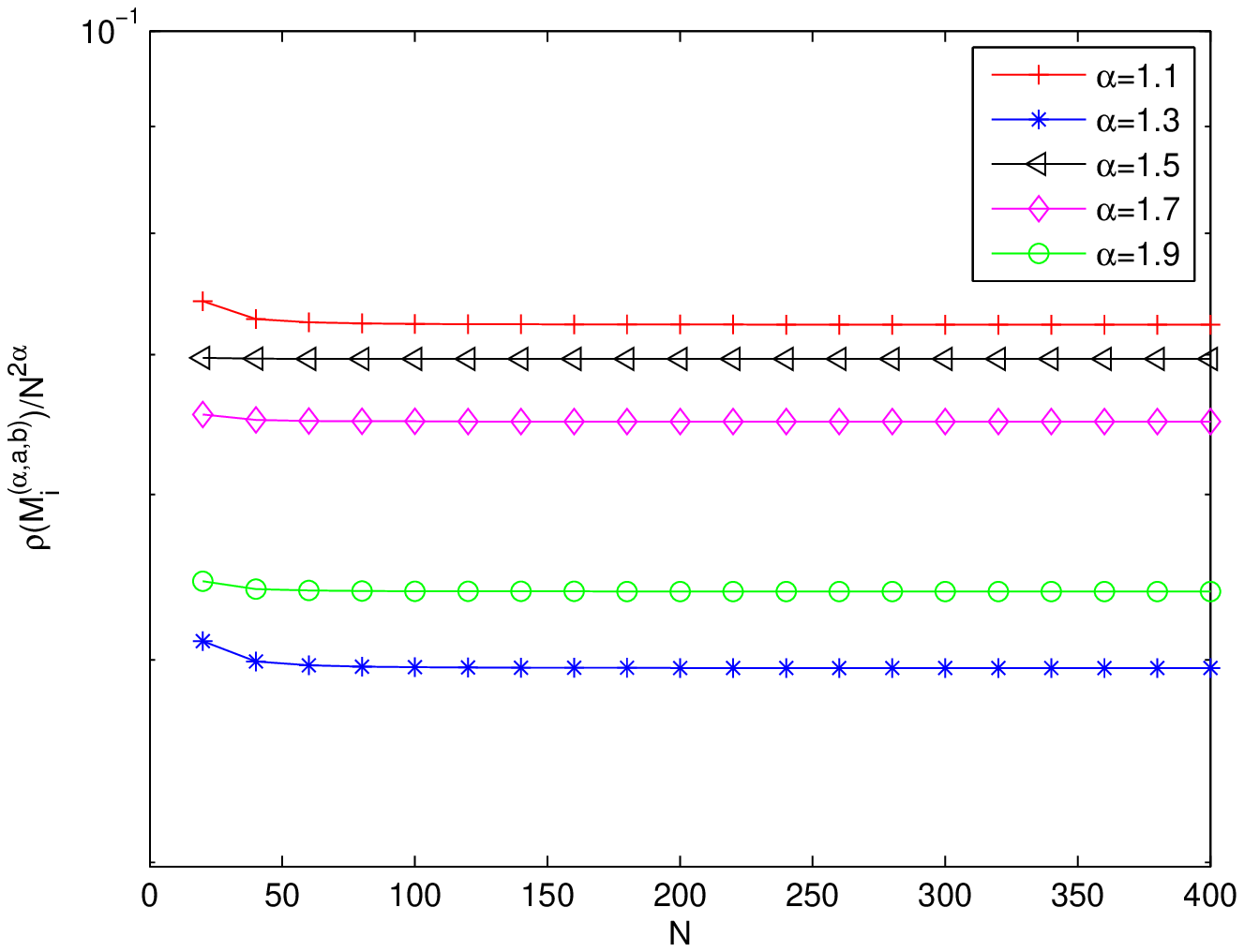}
\par {(c)  $M_{i}^{(\alpha,a,b)}=M_{RL,L}^{(\alpha,-\frac{1}{2},-\frac{1}{2})}$.}
\end{minipage}
\begin{minipage}{0.4\textwidth}\centering
\includegraphics[scale=0.35]{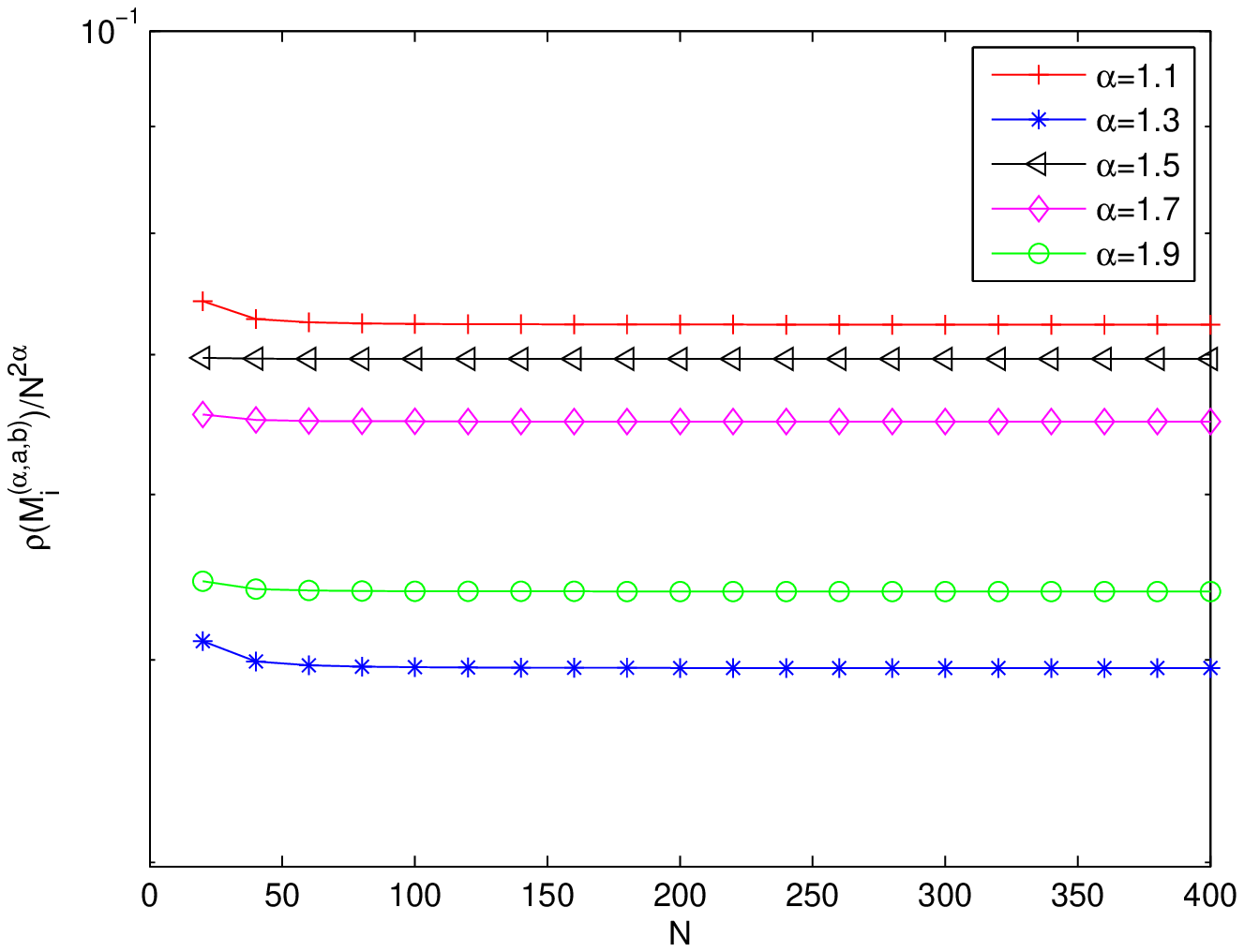}
\par {(d)  $M_{i}^{(\alpha,a,b)}=M_{RL,R}^{(\alpha,-\frac{1}{2},-\frac{1}{2})}$.}
\end{minipage}
\begin{minipage}{0.4\textwidth}\centering
\includegraphics[scale=0.35]{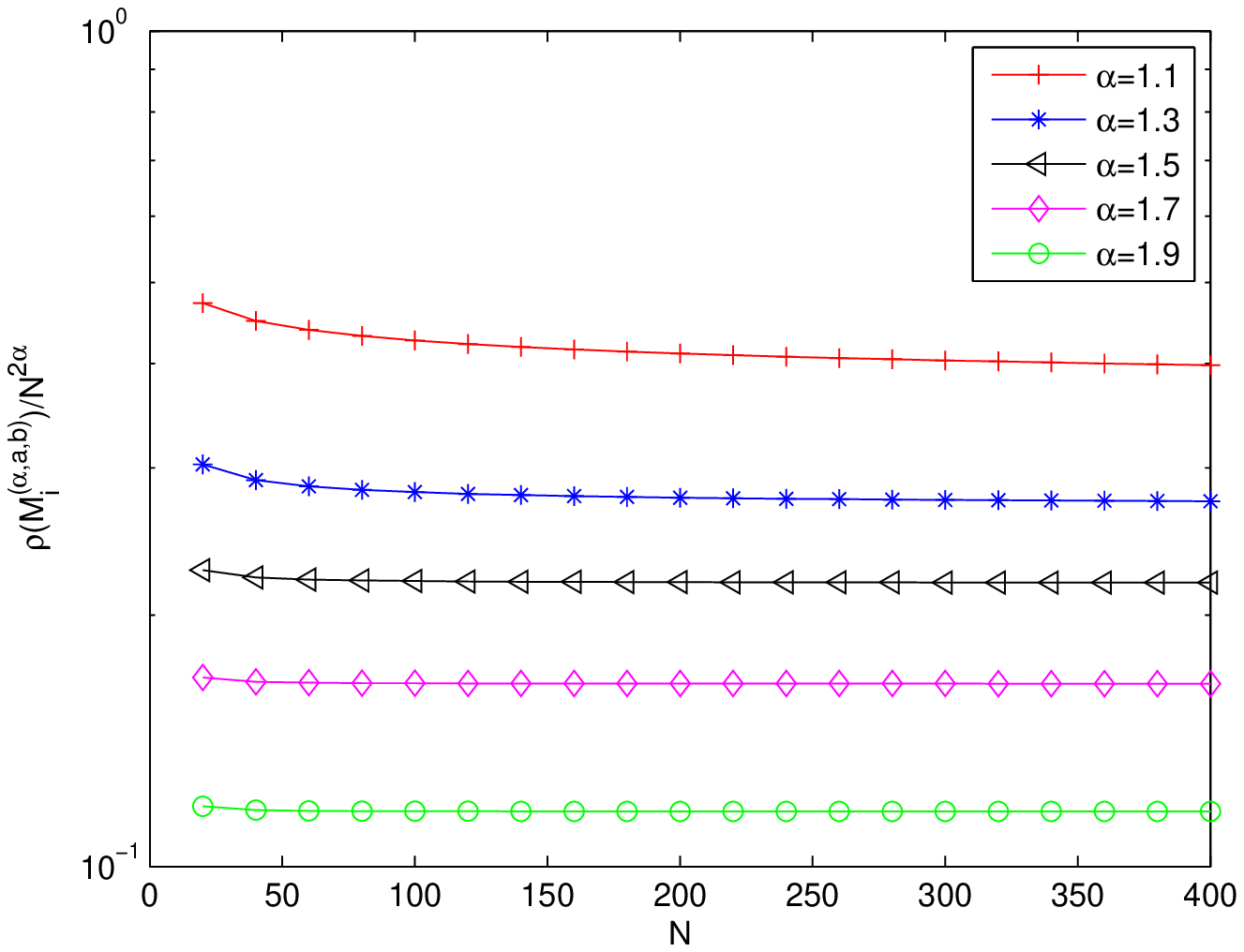}
\par {(e) $M_{i}^{(\alpha,a,b)}=M_{RC}^{(\alpha,-\frac{1}{2},-\frac{1}{2})}$.  }
\end{minipage}
\begin{minipage}{0.4\textwidth}\centering
\includegraphics[scale=0.35]{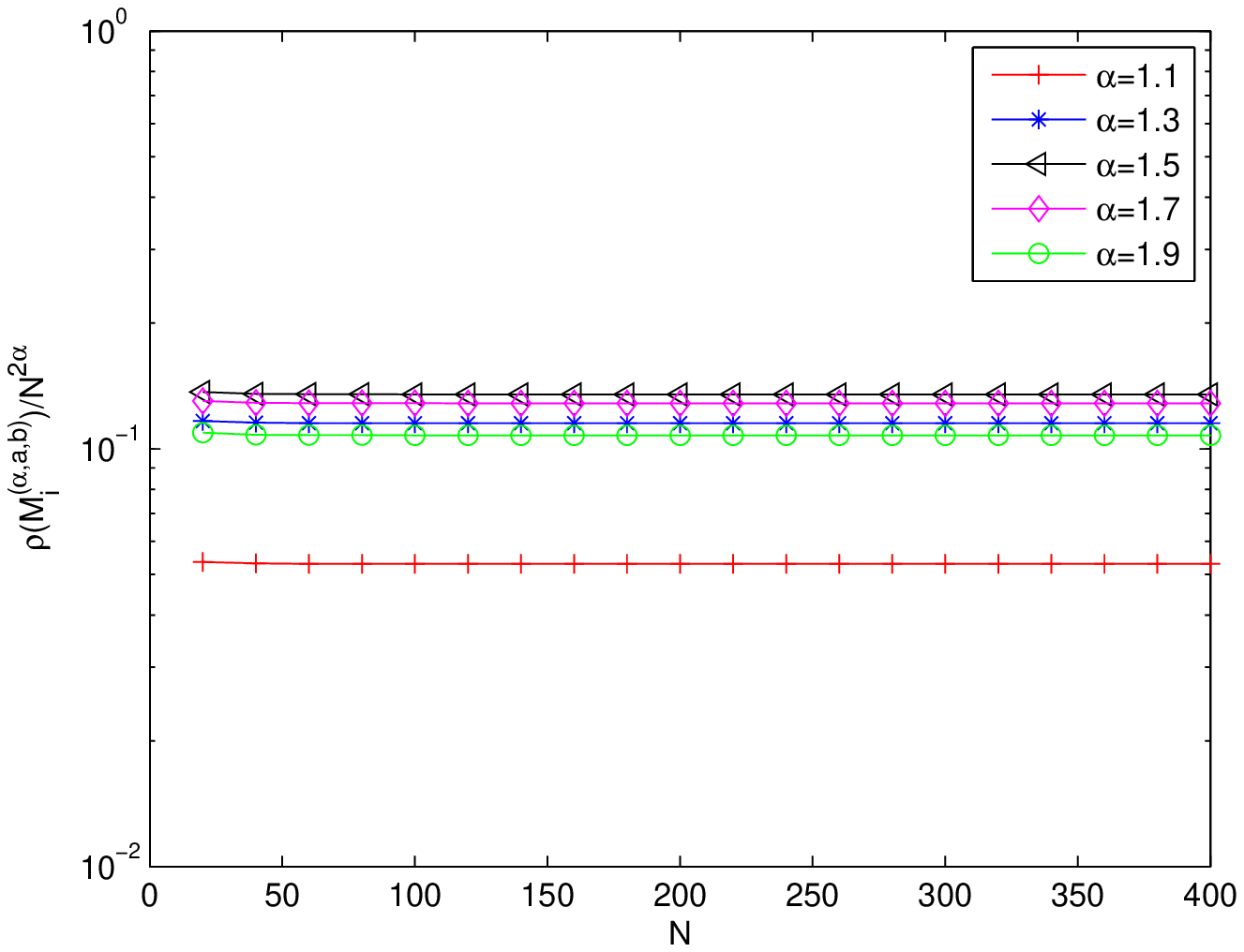}
\par {(f) $M_{i}^{(\alpha,a,b)}=M_{RZ}^{(\alpha,-\frac{1}{2},-\frac{1}{2})}$. }
\end{minipage}
\end{center}
\caption{The boundedness of $\rho\left(M_{i}^{(\alpha,a,b)}\right)/N^{2\alpha}$ for $a=b=-\frac{1}{2}$.\label{fig2}}
\end{figure}

\begin{figure}
\begin{center}
\begin{minipage}{0.4\textwidth}\centering
\includegraphics[scale=0.35]{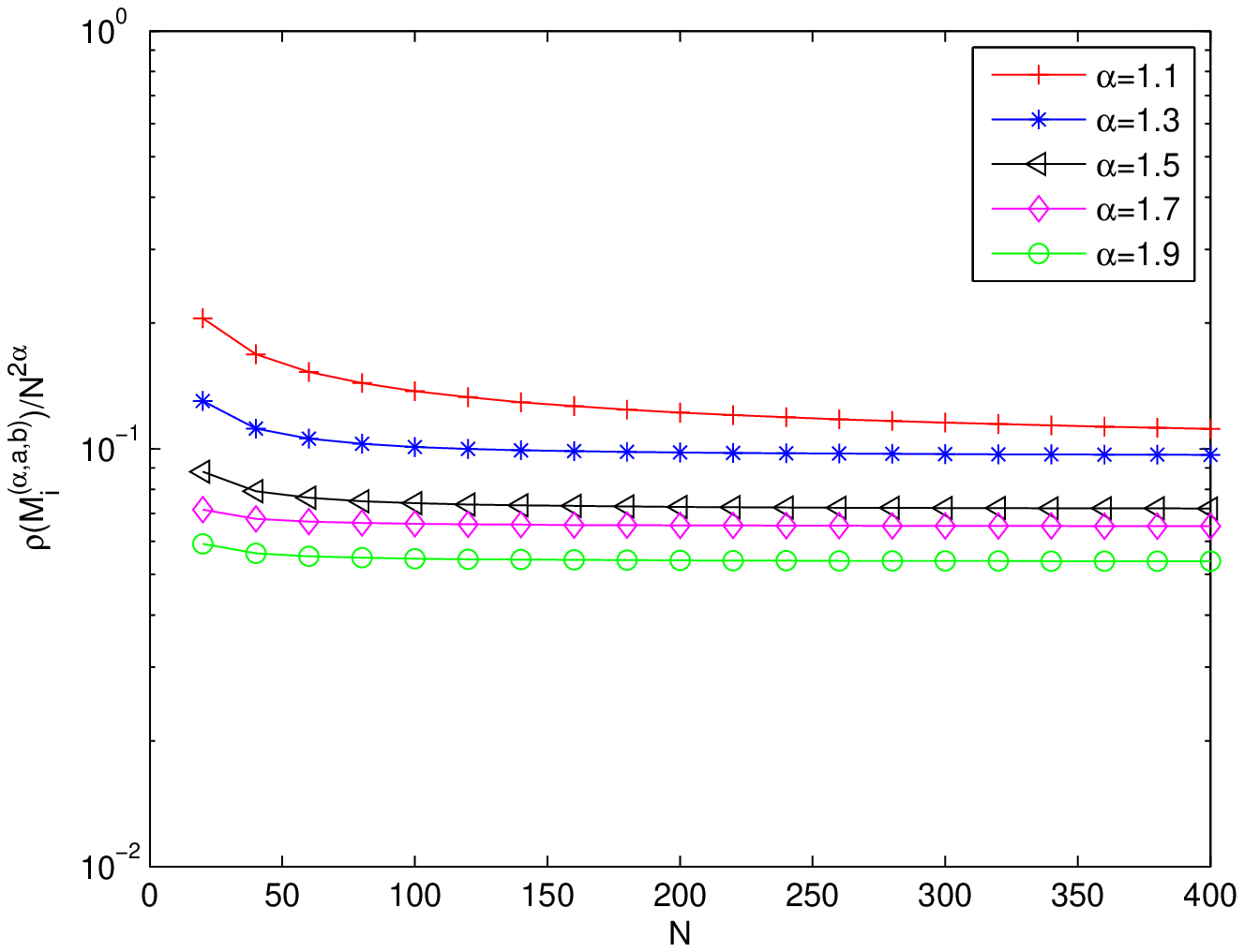}
\par {(a)  $M_{i}^{(\alpha,a,b)}=M_{C,L}^{(\alpha,-\frac{1}{2},\frac{1}{2})}$.}
\end{minipage}
\begin{minipage}{0.4\textwidth}\centering
\includegraphics[scale=0.35]{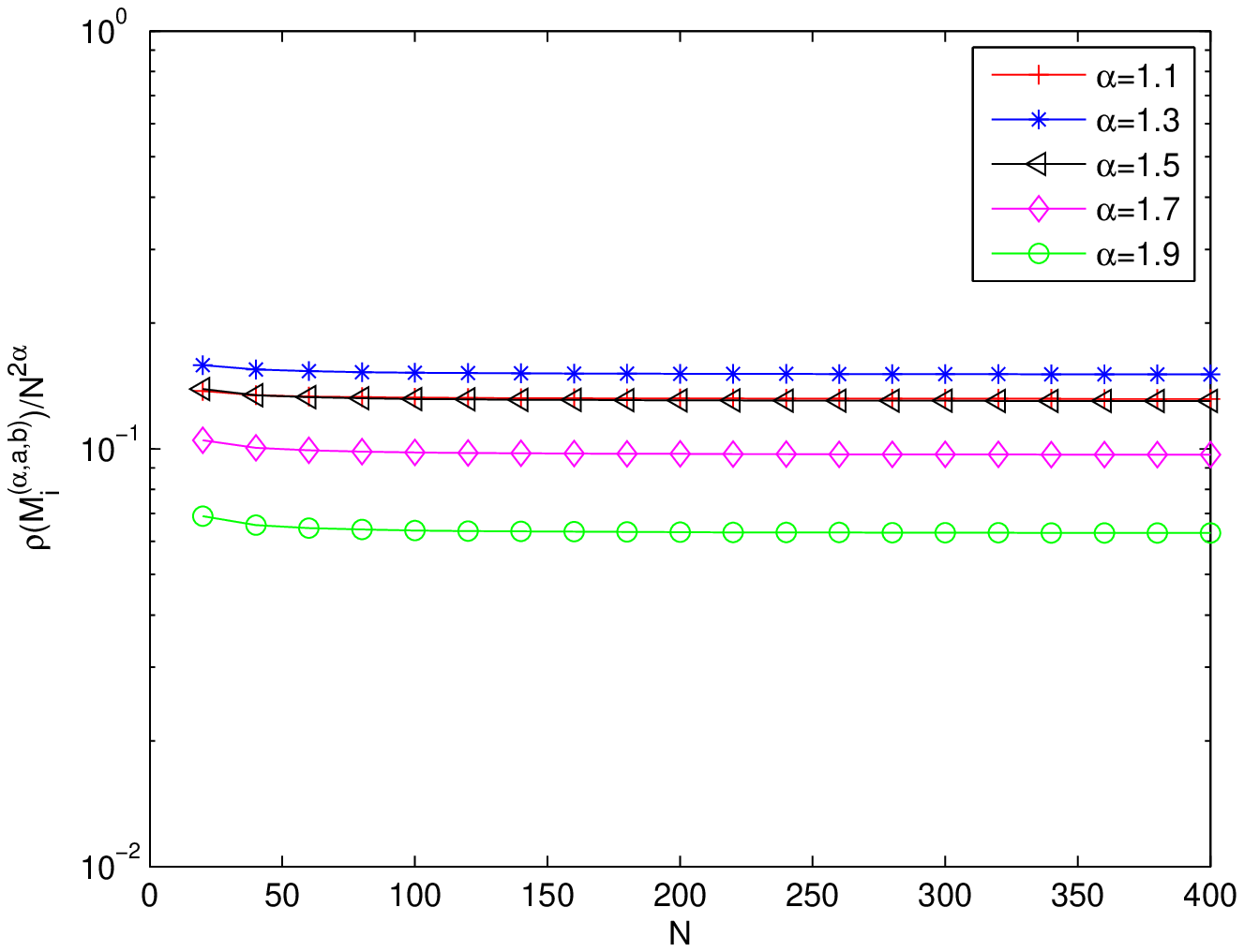}
\par {(b) $M_{i}^{(\alpha,a,b)}=M_{C,R}^{(\alpha,-\frac{1}{2},\frac{1}{2})}$.}
\end{minipage}\\
\begin{minipage}{0.4\textwidth}\centering
\includegraphics[scale=0.35]{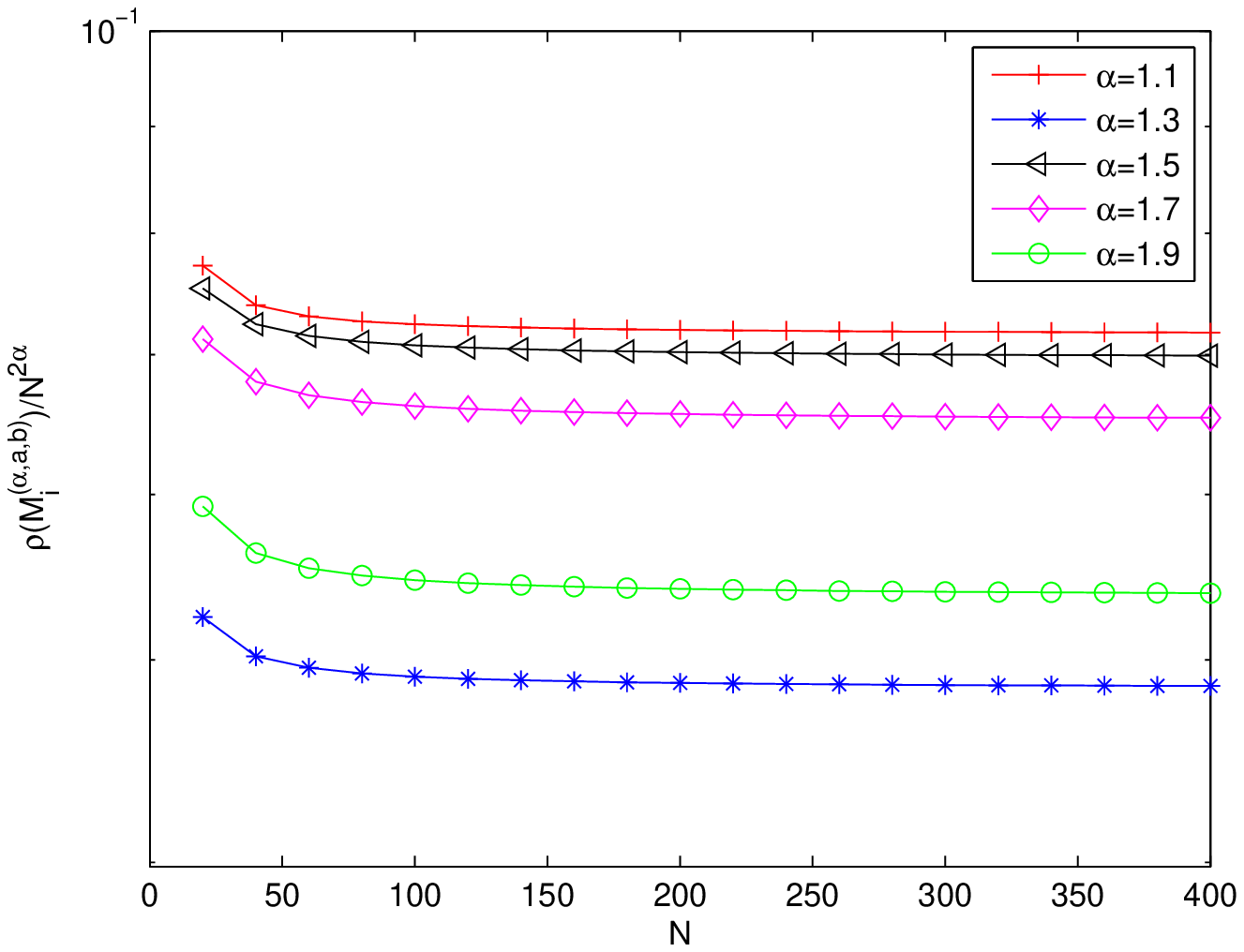}
\par {(c)  $M_{i}^{(\alpha,a,b)}=M_{RL,L}^{(\alpha,-\frac{1}{2},\frac{1}{2})}$.}
\end{minipage}
\begin{minipage}{0.4\textwidth}\centering
\includegraphics[scale=0.35]{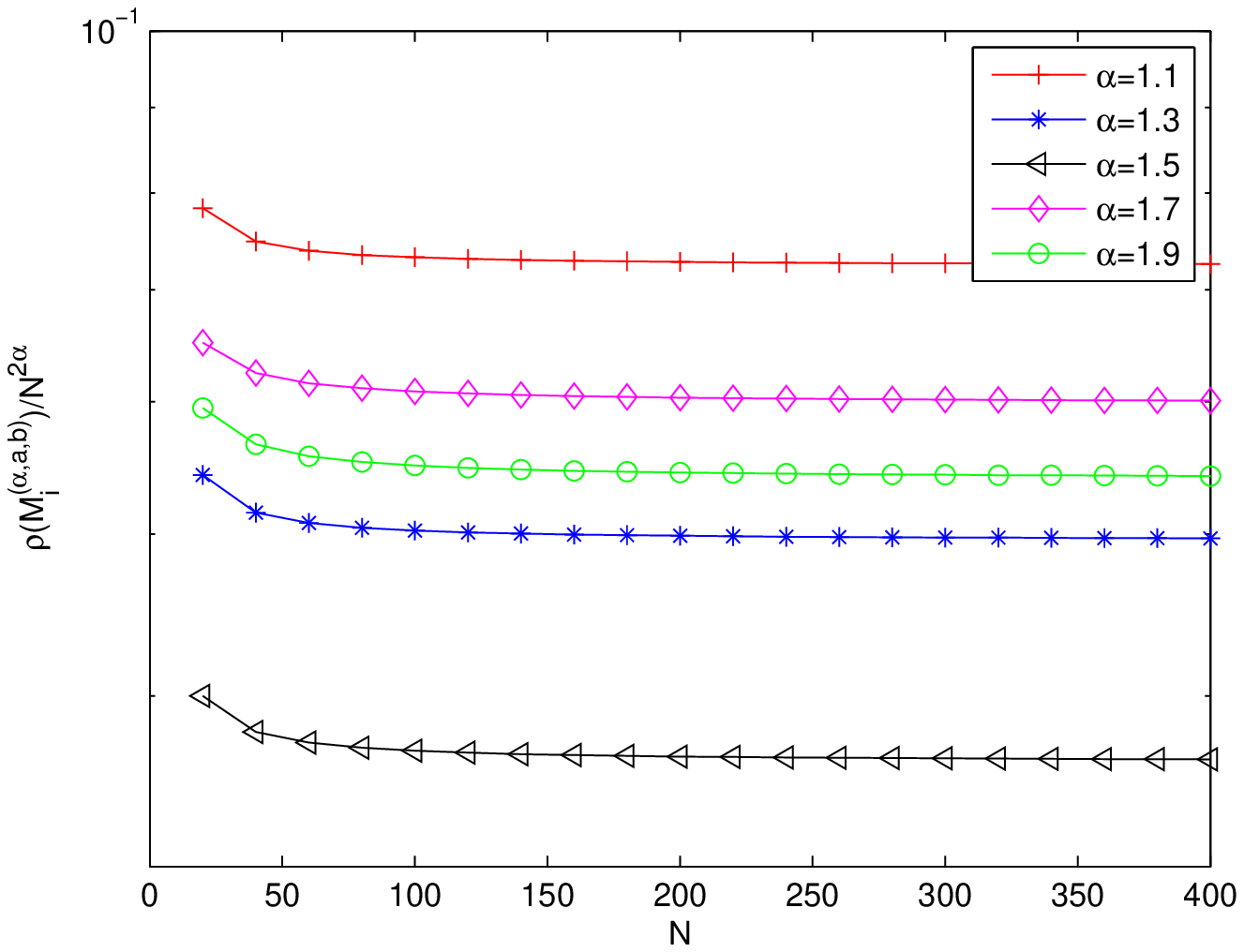}
\par {(d)  $M_{i}^{(\alpha,a,b)}=M_{RL,R}^{(\alpha,-\frac{1}{2},\frac{1}{2})}$.}
\end{minipage}
\begin{minipage}{0.4\textwidth}\centering
\includegraphics[scale=0.35]{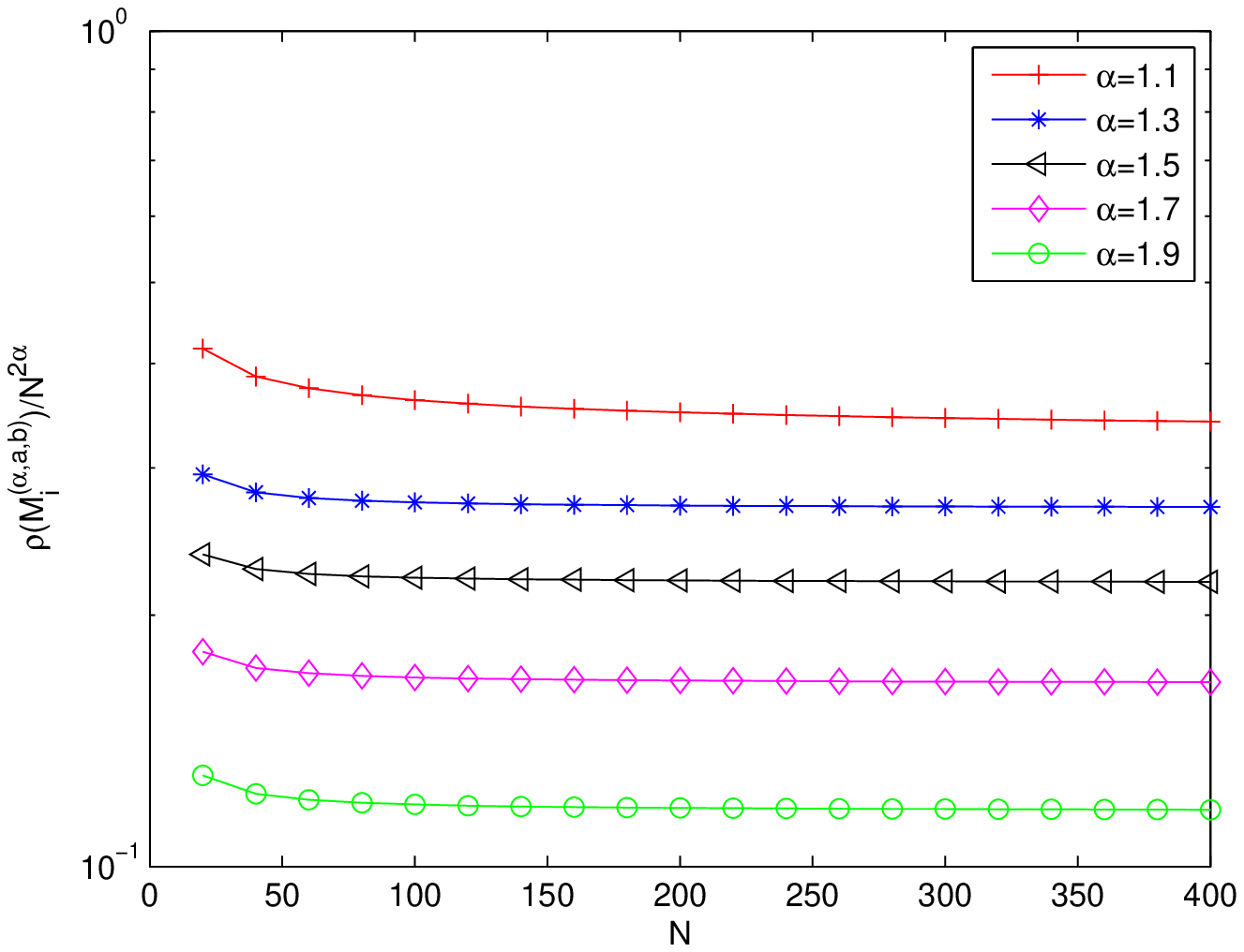}
\par {(e) $M_{i}^{(\alpha,a,b)}=M_{RC}^{(\alpha,-\frac{1}{2},\frac{1}{2})}$.  }
\end{minipage}
\begin{minipage}{0.4\textwidth}\centering
\includegraphics[scale=0.35]{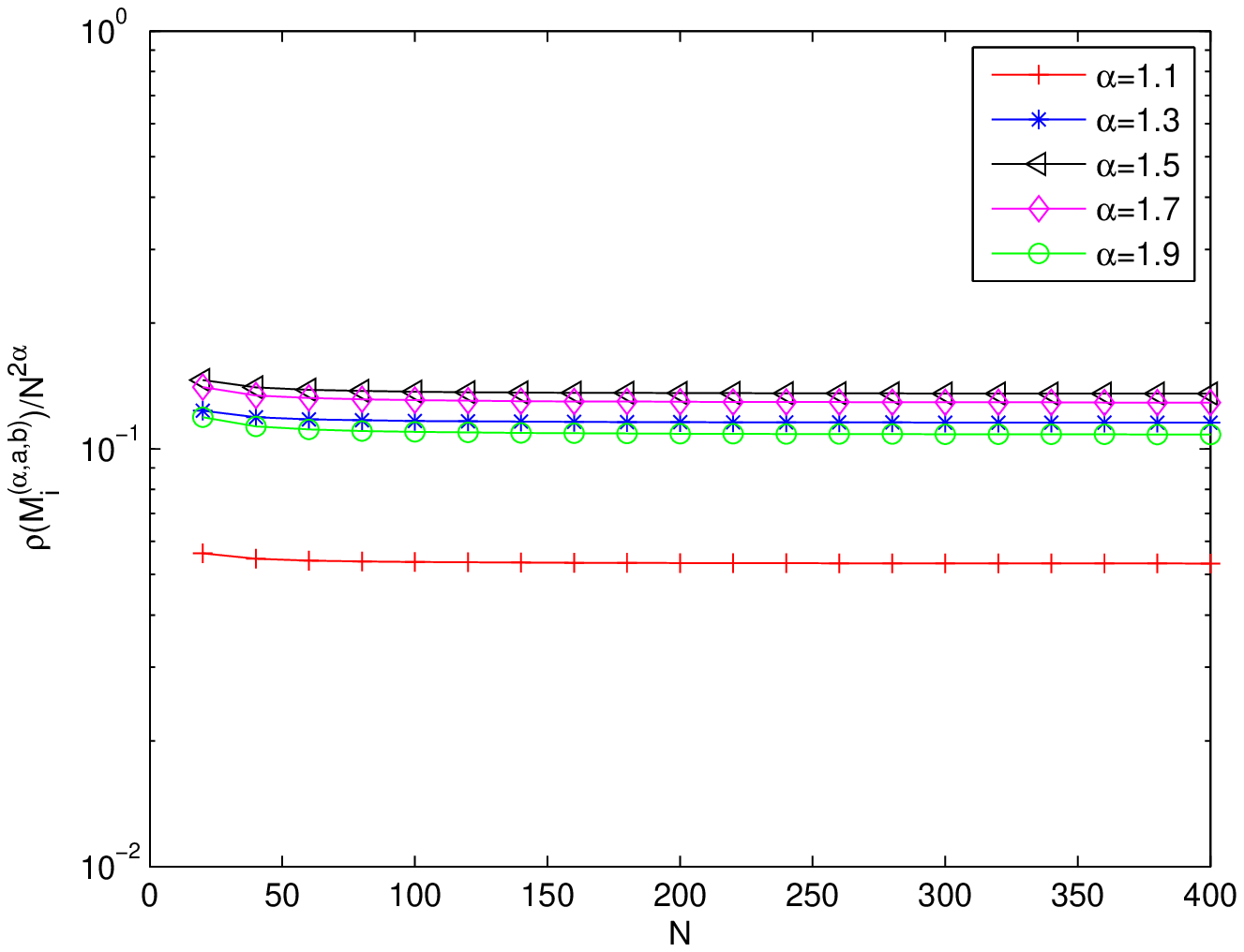}
\par {(f) $M_{i}^{(\alpha,a,b)}=M_{RZ}^{(\alpha,-\frac{1}{2},\frac{1}{2})}$. }
\end{minipage}
\end{center}
\caption{The boundedness of $\rho\left(M_{i}^{(\alpha,a,b)}\right)/N^{2\alpha}$ for $a=-\frac{1}{2},b=\frac{1}{2}$.\label{fig3}}
\end{figure}

In  Figures \ref{fig1}--\ref{fig3},
we plot the behaviors of $\rho\left(M_{i}^{(\alpha,a,b)}\right)/N^{2\alpha}$ with respect to $N$
for different $(a,b)$ ($(a,b)=(0,0),(a,b)=(-1/2,-1/2),(a,b)=(-1/2,1/2)$) and different fractional
order $\alpha\,(\alpha=1.1,1.3,1.5,1.7,1.9)$. Obviously, $\rho\left(M_{i}^{(\alpha,a,b)}\right)/N^{2\alpha}$
is bounded in such  cases. We also test the corresponding  cases of
$0<\alpha<1$, which show similar behaviors.
Here, we conjecture that \eqref{eq:fint20} holds for all $\alpha$.

\begin{remark}
If $\alpha$ is reduced to a positive integer, i.e., $\alpha=m$, then the matrices
$M_{i}^{(\alpha,a,b)}$ are reduced the classical differential matrix $D^{(m,a,b)}$.
\end{remark}

\begin{remark}
If the Jacobi-Gauss or Jacobi-Gauss-Radau collocation points are chosen, then we can
similarly derive the corresponding differential matrices as ${}_{C}{{D}}^{(a,b,\alpha)}_L$,
${}_{C}{{D}}^{(a,b,\alpha)}_R$, ${}_{RL}{{D}}^{(a,b,\alpha)}_L$, and ${}_{RL}{{D}}^{(a,b,\alpha)}_R$.
One just need to replace $c_{k,j}$ defined by \eqref{eq:ckj} with
\begin{equation*} \label{eq:ckj-2}
c_{k,j}=\frac{P_k^{a,b}(x_j)\omega_j}{\gamma_k^{a,b}},{\qquad}k=0,1,...,N
\end{equation*}
to obtain the corresponding results \cite{ShenTangWang2011}.
\end{remark}

\section{Applications}
In this section, we illustrate how to use the fractional differential matrices developed
in the previous section to solve fractional differential equations.

\textbf{Application to the fractional ordinary differential equation:}
Consider the following Baglay--Torvic equation
\begin{equation}\label{eq:eg1-1}
u{''}(x)+b(x){}_CD^{\alpha}_{x_a,x}u(x)+c(x)u(x)=f(x),\quad x\in(x_a,x_b),\quad\alpha\in (1,2)
\end{equation}
with the conditions
\begin{itemize}
  \item[(I)] $u(x_a)=\varphi_a, {\quad}u'(x_a)=\varphi'_a$;
  \item[(II)] $u(x_a)=\varphi_a, {\quad}u(x_b)=\varphi_b$.
\end{itemize}

Suppose that $\hat{x}_j\,(j=0,1,...,x_N)$ are JGL points on the interval $[-1,1]$. Then
the corresponding JGL points $x_j\in [x_a,x_b]$ can be obtained as
$$x_j=\frac{(x_b-x_a)\hat{x}_j+x_a+x_b}{2}.$$

Assume that $u_N=u_N(x)=\sum_{j=0}^N\hat{u}_jF_j(x)$ is the approximation of $u(x),x\in[x_a,x_b]$.
Then we replace $u(x)$
in \eqref{eq:eg1-1} with $u_N(x)$ and let $x=x_j\,(j=1,2,...,N-1)$, which yields
\begin{equation}\label{eq:eg1-2}
(u_N(x)){''}(x_j)+\left[{}_CD^{\alpha}_{x_a,x}u_N(x)\right]_{x=x_j}+u_N(x_j)=f(x_j).
\end{equation}
Note that
\begin{equation}\label{eq:eg1-3}
{}_CD^{\alpha}_{x_a,x}u_N(x)
=\left(\frac{x_b-x_a}{2}\right)^{-\alpha}{}_CD^{\alpha}_{-1,\hat{x}}u_N(\hat{x}).
\end{equation}
Hence we can obtain from \eqref{eq:fint12}, \eqref{eq:fracmat}, and \eqref{eq:eg1-3}
\begin{equation}\label{eq:eg1-4}
\left(\begin{array}{c}
  {}_{C}D_{-1,\hat{x}_0}^{\alpha}u_N(\hat{x}_0) \\
  {}_{C}D_{-1,\hat{x}_1}^{\alpha}u_N(\hat{x}_1) \\
    \vdots\\
 {}_{C}D_{-1,\hat{x}_{N}}^{\alpha}u_N(\hat{x}_{N})
\end{array}\right)=\left({}_C{D}^{(a,b,\alpha)}_{L,x_a,x_b}\right)(\hat{u}_0,\hat{u}_1,\cdots,\hat{u}_N)^T,
\end{equation}
where ${}_C{D}^{(a,b,\alpha)}_{L,x_a,x_b} \in \mathbb{R}^{(N+1)\times(N+1)}$  satisfying
\begin{equation}\label{eq:eg1-5}
{}_C{D}^{(a,b,\alpha)}_{L,x_a,x_b} =\left(\frac{x_b-x_a}{2}\right)^{-\alpha}\left({}_C{D}^{(a,b,\alpha)}_L\right),
\end{equation}
with ${}_C{D}^{(a,b,\alpha)}_L$ defined by \eqref{eq:fracmat}.

Combing the initial condition (I) in \eqref{eq:eg1-1} and \eqref{eq:eg1-4}, we derive the
following linear system
\begin{equation}\label{eq:eg1-6}
A^{(\alpha,a,b)}\mathbf{\hat{u}}=\mathbf{f},
\end{equation}
where $\mathbf{f}=(f({x}_1),f({x}_2),...,f({x}_{N-1}),\varphi_a')^T$,
$\mathbf{\hat{u}}=(\hat{u}_0,\hat{u}_1,...,\hat{u}_N)$, $\hat{u}_0=u(0)=\varphi_a$ is known,
and $A^{(\alpha,a,b)} \in \mathbb{R}^{N\times(N+1)}$  satisfying
\begin{equation*}\label{eq:eg1-7}
\begin{aligned}
&\left(A^{(\alpha,a,b)}\right)_{i-1,j}
=\left({}_C{D}^{(a,b,2)}_{L,x_a,x_b}\right)_{i,j}
+b(x_{i})\left({}_C{D}^{(a,b,\alpha)}_{L,x_a,x_b}\right)_{i,j}+c(x_{i})(E)_{i,j},
{\quad}i=1,2,...,N-1,j=0,1,...N,\\
&\left(A^{(\alpha,a,b)}\right)_{N-1,j}=\left({}_C{D}^{(a,b,1)}_{L,x_a,x_b}\right)_{0,j}
=\left(\frac{x_b-x_a}{2}\right)^{-1}\left({D}^{(a,b,1)}\right)_{0,j},
{\quad}j=0,1,...,N,
\end{aligned}
\end{equation*}
in which $E$ is an ${(N+1)\times(N+1)}$  identity matrix and ${D}^{(a,b,1)}$ is the first-order
differential matrix.

If we use the boundary conditions (II) in \eqref{eq:eg1-1}, we can obtain the following
linear system
\begin{equation}\label{eq:eg1-8}
A^{(\alpha,a,b)}\mathbf{\hat{u}}=\mathbf{f},
\end{equation}
where $\mathbf{f}=(f({x}_1),f({x}_2),...,f({x}_{N-1}))^T$,
$\mathbf{\hat{u}}=(\hat{u}_0,\hat{u}_1,...,\hat{u}_N)$, $\hat{u}_0=u(x_a)=\varphi_a$ and
$\hat{u}_N=u(x_b)=\varphi_b$ are known,
and $A^{(\alpha,a,b)} \in \mathbb{R}^{(N-1)\times(N+1)}$  satisfying
\begin{equation*}\label{eq:eg1-9}
\begin{aligned}
&\left(A^{(\alpha,a,b)}\right)_{i-1,j}
=\left({}_C{D}^{(a,b,2)}_{L,x_a,x_b}\right)_{i,j}
+b(x_{i})\left({}_C{D}^{(a,b,\alpha)}_{L,x_a,x_b}\right)_{i,j}+c(x_{i})(E)_{i,j},
{\quad}i=1,2,...,N-1,j=0,1,...,N.
\end{aligned}
\end{equation*}

\textbf{Application to the fractional partial differential equations}:
Consider the following fractional diffusion equation
\begin{equation}\label{eq:eg2-1}
\left\{\begin{aligned}
&\frac{\partial u}{\partial t}=c_{+}(x){}_{RL}D_{x_a,x}^{\alpha}u(x,t)+
c_{-}(x){}_{RL}D_{x,x_b}^{\alpha}u(x,t)+f(x,t),
{\quad}(x,t){\,\in\,}(x_a,x_b){\times}(0,T],T>0,\\
&u(x,0)=\phi_0(x),{\quad}x{\,\in\,}[x_a,x_b],\\
&u(x_a,t)=\varphi_a(t),{\quad}u(x_b,t)=\varphi_b(t){\quad}t{\,\in\,}\partial{\Omega}\times(0,T],
\end{aligned}\right.
\end{equation}

Assume that  $u_N(x,t)$ be an $N$th-order polynomial with respect $x$  for fixed $t$,
$x_j\,(j=0,1,...,N)$ are the JGL collocation points on $[x_a,x_b]$.
Inserting $u_N(x,t)$ into  \eqref{eq:eg2-1} and letting $x=x_j\,(j=1,2,..,N-1)$ yield
\begin{equation}\label{eq:eg2-2}
\frac{\mathrm{d} u_N(x_j,t)}{\dx[t]}=c_{+}(x_j){}_{RL}D_{x_a,x_j}^{\alpha}u_N(x_j,t)+
c_{-}(x_j){}_{RL}D_{x_j,x_b}^{\alpha}u_N(x_j,t)+f(x_j,t).
\end{equation}
As is done for the fractional ordinary differential equations, we can obtain the
matrix representation of \eqref{eq:eg2-2} below
\begin{equation}\label{eq:eg2-3}
\frac{\mathrm{d} \mathbf{\hat{u}}(t)}{\dx[t]}=A^{(\alpha,a,b)}\mathbf{\hat{u}}(t)+\mathbf{f}(t).
\end{equation}
where $\mathbf{\hat{u}}(t)=(u_N(x_0,t),u_N(x_1,t),...,u_N(x_N,t))^T$,
$u_N(x_0,t)=\varphi_a(t)$ and
$u_N(x_N,t)=\varphi_b(t)$ are known functions,
$\mathbf{f}(t)=(f(x_1,t),f(x_2,t),...,f(x_{N-1},t))^T$,
and
$A^{(\alpha,a,b)} \in \mathbb{R}^{(N-1)\times(N+1)}$  satisfying
\begin{equation*}\label{eq:eg2-4}
\begin{aligned}
\left(A^{(\alpha,a,b)}\right)_{i-1,j}
=&\left(\frac{x_b-x_a}{2}\right)^{-\alpha}
\left[c_{+}(x_{i})\left({}_{RL}{D}^{(a,b,\alpha)}_{L}\right)_{i-1,j}
+c_{-}(x_{i})\left({}_{RL}{D}^{(a,b,\alpha)}_{R}\right)_{i,j}\right],\\
&i=1,2,...,N-1,j=0,1,...,N,
\end{aligned}
\end{equation*}
in which  the matrices ${}_{RL}{D}^{(a,b,\alpha)}_{L}$ and  ${}_{RL}{D}^{(a,b,\alpha)}_{R}$ are
defined  by \eqref{eq:fint14} and \eqref{eq:fint16}, respectively.
From \eqref{eq:eg2-1}, the initial condition of \eqref{eq:eg2-3} can be given as follows
\begin{equation}\label{eq:eg2-5}
\mathbf{\hat{u}}(0)=(\phi_0(x_0),\phi_0(x_1),...,\phi_0(x_{N}))^T.
\end{equation}

Now, the fractional ordinary system \eqref{eq:eg2-3} with initial condition
\eqref{eq:eg2-5} can be solved by any known methods as the Euler method, the trapezoidal
rule, linear multi-step methods, Runge-Kutta methods, and so on.
In the numerical simulation, we use the trapezoidal rule to solve \eqref{eq:eg2-3}.


\section{Numerical examples}
This section provides the numerical examples to verify the methods
obtained in the preceding sections.

\begin{example}\label{eg1}
Consider the following Baglay-Torvik equation \cite{DohaBhrawyEzz-Eldien2011,Podlubny1999}
\begin{equation}\label{eqeg1}
u{''}(x)+{}_CD^{\alpha}_{0,x}u(x)+u(x)=f(x),\quad
x\in(0,1],{\quad}1<\alpha<2,
\end{equation}
with the initial conditions
\begin{equation}\label{eqeg3.1}
u(0)=0,\quad u'(0)=w.
\end{equation}
Choosing appropriate $f(x)$ such that \eqref{eqeg1} has the exact
solution $u(x)=\sin{wx}$.
\end{example}

We  set $a=b=0$ (Legendre collocation method) and $a=b=-1/2$ (Chebyshev collocation method),
the results are shown in Tables \ref{tb1-1} and \ref{tb1-2},  respectively,
where the maximum absolute errors of the
present method \eqref{eq:eg1-6}  and the shifted Chebyshev tau (SCT) method
developed in \cite{DohaBhrawyEzz-Eldien2011} are shown.
From this example, we can see that our method shows
more accurate results than the SCT method developed in
\cite{DohaBhrawyEzz-Eldien2011}.

If the initial condition $u'(0)=w$ in \eqref{eqeg3.1} is replaced by
$u(1)=\sin{w}$, we can obtain the boundary value problem. In such a case,
we choose the suitable right hand side function $f(x)$ such that \eqref{eqeg1}
has the exact solution  $u(x)=\sin{wx}$. We use the method \eqref{eq:eg1-8}
to solve the this problem with $w=4\pi,\alpha=1.1,1.25,1.4,1.6,1.75,1.9$.
Tables \ref{tb1-3} and \ref{tb1-4} display  the maximum
absolute errors at Legendre-Gauss-Lobatto ($a=b=0$) points
and Chebyshev-Gauss-Lobatto ($a=b=-\frac{1}{2}$) points, respectively.
Clearly, the satisfactory numerical results are obtained.

\begin{table}[t]\caption{\small The absolute errors for
 Example \ref{eg1} with  different $N$ and $a=b=0$, $\alpha=1.5$.}\label{tb1-1}
\begin{center}
\begin{tabular}{ccccccc}
\hline
  $N$& $w$ &Method \eqref{eq:eg1-6}& SCT \cite{DohaBhrawyEzz-Eldien2011}
  & $w$ &Method \eqref{eq:eg1-6} & SCT \cite{DohaBhrawyEzz-Eldien2011} \\
  \hline
4  &1& 2.4175e-04 & 3.4e-04 & $4\pi$ & 1.2516e+01 & 3.9e-00  \\
8  &1& 7.3967e-10 & 4.3e-07 & $4\pi$ & 1.3775e+00 & 4.7e-01  \\
16 &1& 6.4948e-14 & 1.8e-08 & $4\pi$ & 8.5461e-05 & 3.5e-05  \\
32 &1& 2.3959e-13 & 7.1e-10 & $4\pi$ & 4.5841e-12 & 1.4e-06  \\
48 &1& 3.0032e-13 & 9.9e-11 & $4\pi$ & 2.3109e-12 & 1.9e-07  \\
64 &1& 1.5175e-12 & 2.4e-11 & $4\pi$ & 1.0300e-11 & 4.8e-08  \\
 \hline
\end{tabular}
\end{center}
\end{table}

\begin{table}[t]\caption{\small The absolute errors for
 Example \ref{eg1} with  different $N$ and $a=b=-1/2$, $\alpha=1.5$.}\label{tb1-2}
\begin{center}
\begin{tabular}{ccccccc}
\hline
  $N$& $w$ &Method \eqref{eq:eg1-6}& SCT \cite{DohaBhrawyEzz-Eldien2011}
  & $w$ &Method \eqref{eq:eg1-6} & SCT \cite{DohaBhrawyEzz-Eldien2011} \\
  \hline
4  &1& 1.3548e-04 & 3.4e-04 & $4\pi$ & 1.1894e+01 & 3.9e-00  \\
8  &1& 2.8661e-10 & 4.3e-07 & $4\pi$ & 4.8230e-01 & 4.7e-01  \\
16 &1& 5.3291e-15 & 1.8e-08 & $4\pi$ & 2.3177e-05 & 3.5e-05  \\
32 &1& 2.7367e-13 & 7.1e-10 & $4\pi$ & 1.2018e-12 & 1.4e-06  \\
48 &1& 2.4891e-13 & 9.9e-11 & $4\pi$ & 4.5565e-12 & 1.9e-07  \\
64 &1& 4.4387e-13 & 2.4e-11 & $4\pi$ & 7.8381e-12 & 4.8e-08  \\
 \hline
\end{tabular}
\end{center}
\end{table}

\begin{table}[t]\caption{\small The absolute errors for
 Example \ref{eg1} with  boundary value conditions and $a=b=0$, $w=4\pi$.}\label{tb1-3}
\begin{center}
\begin{tabular}{ccccccc}
\hline
  $N$ &$\alpha=1.1$& $\alpha=1.25$ & $\alpha=1.4$ & $\alpha=1.6$ & $\alpha=1.75$ & $\alpha=1.9$    \\
  \hline
8  &2.4934e-01 &2.1828e-01&1.8398e-01 &1.4161e-01&1.1189e-01&7.0096e-02\\
12 &4.4902e-03 &3.8408e-03&3.1023e-03 &2.1753e-03&1.5838e-03&8.9844e-04\\
16 &1.5053e-05 &1.2813e-05&1.0230e-05 &6.9318e-06&4.8828e-06&2.7329e-06\\
20 &1.6351e-08 &1.3898e-08&1.1044e-08 &7.3458e-09&5.0633e-09&2.8151e-09\\
24 &7.5480e-12 &6.2294e-12&5.1675e-12 &3.1971e-12&1.4853e-12&2.6823e-12\\
28 &2.7034e-13 &2.5047e-13&4.1245e-13 &8.9195e-13&4.8161e-13&3.4023e-13\\
 \hline
\end{tabular}
\end{center}
\end{table}

\begin{table}[t]\caption{\small The absolute errors for
 Example \ref{eg1} with  boundary value conditions and $a=b=-\frac{1}{2}$, $w=4\pi$.}\label{tb1-4}
\begin{center}
\begin{tabular}{ccccccc}
\hline
  $N$ &$\alpha=1.1$& $\alpha=1.25$ & $\alpha=1.4$ & $\alpha=1.6$ & $\alpha=1.75$ & $\alpha=1.9$    \\
  \hline
8  & 1.2224e-01 &1.1248e-01&1.0122e-01 &8.6472e-02&7.5269e-02 & 5.8414e-02\\
12 & 1.4914e-03 &1.3028e-03&1.0830e-03 &7.9562e-04&6.0160e-04 & 3.7026e-04\\
16 & 4.2450e-06 &3.6549e-06&2.9657e-06 &2.0617e-06&1.4783e-06 & 8.5876e-07\\
20 & 4.1101e-09 &3.5201e-09&2.8283e-09 &1.9126e-09&1.3308e-09 & 7.6665e-10\\
24 & 1.6886e-12 &1.2088e-12&1.4165e-12 &1.3994e-12&5.8165e-13 & 3.7370e-13\\
28 & 2.3892e-13 &8.9373e-14&6.9333e-13 &3.6260e-13&5.0826e-13 & 1.1949e-12\\
 \hline
\end{tabular}
\end{center}
\end{table}

\begin{example}\label{eqeg2}
Consider the following diffusion equation \cite{CelikDuman2012}
\begin{equation}\label{eg2}
\left\{\begin{aligned}
&\frac{\partial u}{\partial t}={}_{RZ}D_{x}^{\alpha}u(x,t)+f(x,t),
{\quad}(x,t){\,\in\,}(0,1){\times}(0,T],T>0,\\
&u(x,0)=x^2(1-x^2),{\quad}x{\,\in\,}[0,1],\\
&u(0,t)=u(1,t)=0{\quad}t{\,\in\,}\partial{\Omega}\times(0,T].
\end{aligned}\right.
\end{equation}
Choosing the suitable right hand side function $f(x,t)$ such that \eqref{eg2}
has the following exact solution
$$u(x,t)=(t+1)^{\alpha}x^2(1-x)^2.$$
\end{example}

Since the exact solution is a polynomial of order 4, so we choose $N=4$ in the computation.
The time step size ($\tau=10^{-2}$) is the same as that in \cite{CelikDuman2012},
the fractional orders are chosen as $\alpha=1.1,1.3,1.5,1.7,1.9$, the absolute maximum errors
are shown in Table \ref{tb2-1}. The notation $(a,b)=(0,0)$ means the method \eqref{eq:eg2-3}--\eqref{eq:eg2-5}
with $a=b=0$ is used to solve \eqref{eg2}, which is the same for $(a,b)=(-\frac{1}{2},-\frac{1}{2})$
and $(a,b)=(-\frac{1}{2},\frac{1}{2})$.
$N=100$ and $N=200$ in Table \ref{tb2-1} imply the space step size
$h=1/N$ was used in \cite{CelikDuman2012}. Obviously, the present methods show better numerical
solutions here.
\begin{table}[t]\caption{\small The absolute errors for Example \ref{eqeg2} at $t=10$.}\label{tb2-1}
\begin{center}
\begin{tabular}{lccccc}
\hline
Methods &$\alpha=1.1$ & $\alpha=1.3$ & $\alpha=1.5$  & $\alpha=1.7$ &  $\alpha=1.9$\\
  \hline
$(a,b)=(0,0)$                       &1.0120e-8 & 5.7515e-8 & 1.7774e-7 & 4.5405e-7 & 1.0518e-6\\
$(a,b)=(-\frac{1}{2},-\frac{1}{2})$ &1.0116e-8 & 5.7501e-8 & 1.7772e-7 & 4.5402e-7 & 1.0518e-6\\
$(a,b)=(-\frac{1}{2},\frac{1}{2})$  &9.6527e-9 & 5.4857e-8 & 1.6953e-7 & 4.3308e-7 & 1.0032e-6\\
CN \cite{CelikDuman2012}($N=100$)   &7.5271e-5 & 1.6136e-4 & 3.5410e-4 & 8.0486e-4 & 1.9026e-3\\
CN \cite{CelikDuman2012}($N=200$)   &2.0479e-5 & 3.8070e-5 & 8.3236e-5 & 1.9120e-4 & 4.6505e-4\\
 \hline
\end{tabular}
\end{center}
\end{table}



\section{Conclusion}
In this paper, we derive the fractional differential matrices
with respect to the Jacobi-Gauss points. The spectral radius of the derived fractional
differential matrices are numerically investigated, which show the behaviors as \eqref{eq:fint20}.
We also develop the spectral collocation schemes to solve the fractional
differential equations, which show good performances. Tian and Deng \cite{TianDeng2013}
also developed the fractional differential matrices, but their method will blow up
when performing on the computers with double precision with $N$ being suitably large,
i.e., when $N>35$, the results are not believable, see Remark 1 in \cite{TianDeng2013}.
Clearly, our method is more stable than that in \cite{TianDeng2013}, which is also
tested by the numerical examples in the present paper, where all the numerical results
are computed on the computer with double precision by Matlab.

\end{document}